\documentclass[10pt]{article}
\usepackage{latexsym,amsmath,amssymb,graphics,amscd, empheq}
\usepackage{color}
\usepackage{soul}
\usepackage{stackengine}
\usepackage[breaklinks,colorlinks,
  citecolor=blue,
  urlcolor=blue]{hyperref}

\addtocounter{footnote}{1}
\makeatletter
\@addtoreset{table}{bsection}
\def\thetable{\thesection.\@arabic\c@table}
\def\fps@table{h, t}
\@addtoreset{equation}{section}

\makeatother

\newtheorem{theorem}{Theorem}[section]
\newtheorem{definition}[theorem]{Definition}
\newtheorem{lemma}[theorem]{Lemma}
\newtheorem{remark}[theorem]{Remark}
\newtheorem{proposition}[theorem]{Proposition}

\newtheorem{corollary}[theorem]{Corollary}

\newcommand{\bfi}{\bfseries\itshape}
\newcommand{\vertiii}[1]{{\left\vert\kern-0.25ex\left\vert\kern-0.25ex\left\vert #1 
    \right\vert\kern-0.25ex\right\vert\kern-0.25ex\right\vert}}
\newsavebox{\savepar}

\makeatletter

\usepackage{scalerel,stackengine}
\stackMath
\newcommand\reallywidehat[1]{%
\savestack{\tmpbox}{\stretchto{%
  \scaleto{%
    \scalerel*[\widthof{\ensuremath{#1}}]{\kern-.6pt\bigwedge\kern-.6pt}%
    {\rule[-\textheight/2]{1ex}{\textheight}}
  }{\textheight}%
}{0.5ex}}%
\stackon[1pt]{#1}{\tmpbox}%
}

\begin{document}

\title{\textbf{Memory and forecasting capacities of nonlinear recurrent networks}}
\author{Lukas Gonon$^{1}$, Lyudmila Grigoryeva$^{2}$, and Juan-Pablo Ortega$^{3, 4}$}
\date{}
\maketitle

\begin{abstract}
The notion of memory capacity, originally introduced for echo state and linear networks with independent inputs, is generalized to nonlinear recurrent networks with stationary but dependent inputs. The presence of dependence in the inputs makes natural the introduction of the network forecasting capacity, that measures the possibility of forecasting time series values using network states. Generic bounds for memory and forecasting capacities are formulated in terms of the number of neurons of the {nonlinear recurrent} network and the autocovariance function {or the spectral density} of the input. These bounds generalize well-known estimates in the literature to a dependent inputs setup. Finally, { for the particular case} of linear recurrent networks with independent inputs it is proved that the memory capacity is given by the rank of the associated controllability matrix, {a fact that has been for a long time assumed to be true without proof by the community}.
\end{abstract}

\bigskip


\makeatletter
\addtocounter{footnote}{1} \footnotetext{%
Ludwig-Maximilians-Universit\"at M\"unchen.  Faculty of Mathematics, Informatics and Statistics. Theresienstrasse 39,
80333 Munich. 
Germany. 
{\texttt{gonon@math.lmu.de}}}
\addtocounter{footnote}{1} \footnotetext{%
Department of Mathematics and Statistics. Universit\"at Konstanz. Box 146. D-78457 Konstanz. Germany. {\texttt{Lyudmila.Grigoryeva@uni-konstanz.de} }}
\addtocounter{footnote}{1} \footnotetext{%
Universit\"at Sankt Gallen. Faculty of Mathematics and Statistics. Bodanstrasse 6.
CH-9000 Sankt Gallen. Switzerland. {\texttt{Juan-Pablo.Ortega@unisg.ch}}}
\addtocounter{footnote}{1} \footnotetext{%
Centre National de la Recherche Scientifique (CNRS). France. }
\makeatother

\medskip

\medskip

\medskip

\section{Introduction}

{\bfi  Memory capacities} have been introduced in~\cite{Jaeger:2002} in the context of recurrent neural networks in general and of echo state networks (ESNs) \cite{Matthews:thesis, Matthews1994, Jaeger04} in particular, as a way to quantify the amount of information contained in the states of a  state-space system in relation with past inputs { and  as a measure of the ability of the network to retain the dynamic features of processed signals}.

In the original definition, the memory capacity was defined as the sum of the coefficients of determination of the different linear regressions that use the state of the system at a given time as covariates and the values of the input at a given { lagged time} in the past as dependent variables. This notion has been the subject of much research in the reservoir computing literature~\cite{White2004, Ganguli2008, Hermans2010, dambre2012, esn2014, linearESN, farkas:bosak:2016, Goudarzi2016, Xue2017, Verzelli2019a} where most of the efforts have been concentrated in linear and echo state systems. Analytical expression of the capacity of time-delay reservoirs have been formulated in \cite{GHLO2014_capacity, RC3} and various proposals for optimized reservoir architectures can be obtained by maximizing  {the capacity as a function of reservoir hyperparameters}  \cite{pesquera2012, GHLO2014, ortin2019tackling, ortin2020delay}. Additionally, memory capacities have been extensively compared with other related concepts like Fisher information-based criteria \cite{Tino2013, Livi2016, tino:symmetric}.

All the above-mentioned works consider exclusively  { independent or white noise} input signals and it is, to our knowledge only in \cite{charles2014short} (using autocorrelated inputs via a sparsity model), \cite{RC4pv} (under strong high-order stationarity assumptions), \cite{Charles2017} (using autocorrelated inputs via low-rank multi-input models), and \cite{marzen:capacity} (for linear recurrent networks and inputs coming from a countable hidden Markov model), that the case involving dependent inputs has been treated. Since most signals that one encounters in applications exhibit some sort of temporal dependence, studying this case in detail is of obvious practical importance. Moreover, the presence of dependence makes pertinent considering not only memory capacities, but also the possibility to forecast time series values using network states, that is, {\bfi  forecasting capacities}. This notion has been introduced for the first time in \cite{marzen:capacity} and studied in detail for linear recurrent networks under the hypothesis that the inputs are realizations of a countable hidden Markov model. { That} work shows, in particular, that linear networks optimized for memory capacity do not necessarily have a good forecasting capacity and vice versa.

{This paper contains two main contributions. First, we  extend  the results on memory and forecasting capacities available in the literature exclusively for either linear recurrent or echo state networks and for  uncorrelated/independent inputs {\it to general non-linear systems and to dependent inputs that are only assumed to be stationary}. We show that under particular assumptions our bounds reduce  to those  well-known in the literature.  More specifically, it is known since~\cite{Jaeger:2002} that the memory capacity of an ESN or a linear recurrent network defined using independent inputs (we call this the {\it classical linear case}) is bounded above by  the  number of its output neurons or, equivalently, by the dimensionality of the corresponding state-space representation. We show that these memory capacity bounds in~\cite{Jaeger:2002} immediately follow from our results. Second, in the linear case we reveal new relations between the memory capacity, spectral properties of the connectivity matrix of the network, and the so-called Kalman's characterization of the controllability of a linear system. More explicitly, \cite{Jaeger:2002}  shows that the memory capacity is maximal if and only if Kalman's controllability rank condition \cite{Kalman2010, sontag1991kalman, sontag:book} is satisfied. In this work we make a step further and prove that {\it  in the classical case the memory capacity is given exactly by the rank of the controllability matrix}.

The paper is organized as follows:
\begin{itemize}
\item Section \ref{Recurrent neural networks with stationary inputs} introduces recurrent neural networks with linear readouts in relation with state-space representations. We focus on a large class of state-space systems that satisfy the so-called echo state property (ESP) and which guarantees that they uniquely determine an input/output system (also referred as {\it filter} in this paper). We recall well-known sufficient conditions for this property to hold, the notions of system morphism and isomorphism, and discuss how the non-uniqueness of state-space representations can be handled. We also carefully introduce the stationarity hypotheses that are invoked in the rest of the paper. Finally,  in Proposition \ref{std rc io}  we introduce an important technical result that shows that if we have a state system and an input for which the output process is covariance stationary and the corresponding covariance matrix is non-singular, then an isomorphic system representation exists whose corresponding state process is standardized, that is,  the states have mean zero and covariance matrix equal to the identity. This standardization  leads to systems that are easier to handle in terms of the computation of memory and forecasting capacities, which is profusely exploited later on in the main results of the paper.
\item Section \ref{Memory and forecasting capacity bounds for stationary inputs} contains the first main contribution of the paper.  We first provide the definitions of the memory and forecasting capacities of nonlinear systems with linear readouts in the presence of  stationary inputs and ouputs. Second, Lemma \ref{capacity invariant linear morphism} shows that the memory and forecasting capacities of state-space systems with linear readouts are invariant with respect to linear system morphisms; this is an technical tool late or in Section \ref{The linear case}. Finally, Theorem \ref{capacity bounds} provides bounds for the memory and the forecasting capacities of generic recurrent networks with linear readouts and with stationary inputs in terms of the dimensionality of the corresponding state-space representation and the autocovariance or the spectral density function of the input, which is assumed to be second-order stationary. These bounds reduce to those in~\cite{Jaeger:2002} when the inputs are independent. 
\item Section \ref{The linear case} is exclusively devoted to the linear case. We study separately the cases in which the (time-independent) covariance matrices of the state process are invertible and non-invertible. In the regular case, explicit expressions for the memory and the forecasting capacities can be stated (see Proposition \ref{the linear case with bs}) in terms of the matrix parameters of the network (the so-called connectivity and input matrices) and the autocorrelation properties of the input. Moreover, in the {\it classical}  case {(linear recurrent network with independent inputs)}, these expressions yield interesting relations (see Proposition \ref{capacities and kalman}) between maximum memory capacity, spectral properties of the connectivity matrix of the network, and the so-called Kalman's characterization of the controllability of a linear state-space system \cite{Kalman2010}. This last condition  has been already mentioned in \cite{Jaeger:2002} in relation with maximal capacity. When the state covariance matrix is singular, a completely different strategy is adopted based on using the invariance of capacities under linear system morphisms. Theorem \ref{expression capacity linear system} proves that {\it the memory capacity of a linear recurrent network with independent inputs is given by the rank of its controllability matrix}. This statement obviously generalizes the one established in \cite{Jaeger:2002}. Even though, to our knowledge, this is the first rigorous proof of the relation between network memory and the rank of the controllability matrix, that link has been for a long time  part of the reservoir computing folklore. In particular, recent contributions are dedicated to the design of ingenuous configurations that maximize that rank \cite{Rodan2011, Aceituno2017, Verzelli2020}. 
\item Section \ref{Conclusions} concludes the paper and all the proofs are contained in the Appendices in Section \ref{Appendices}.
\end{itemize}}

\section{Recurrent neural networks with stationary inputs}
\label{Recurrent neural networks with stationary inputs}

The results in this paper apply to recurrent neural networks determined by state-space equations of the form:
\begin{empheq}[left={\empheqlbrace}]{align}
\mathbf{x} _t &=F(\mathbf{x} _{t-1}, {\bf z}_t),\label{rc state eq}\\
{\bf  y} _t &= h( \mathbf{x} _t) := \mathbf{W}^\top\mathbf{x} _t + {\bf a}, \label{rc readout eq}
\end{empheq}
for any $t \in \Bbb Z $. These two relations form a {\bfi  state-space system}, where the map $F:  D_N \subset \mathbb{R} ^N \times D _d \subset \mathbb{R}^d \longrightarrow D_N \subset \mathbb{R}^N$, $N, d \in \mathbb{N} $, is called the {\bfi  state map} and  $h: \mathbb{R}^N \longrightarrow \mathbb{R}^m$ the {\bfi  readout } or {\bfi  observation } map that, all along this paper, will be assumed to be affine, that is, it is determined just by a matrix $\mathbf{W} \in \mathbb{M}_{N,m} $ and a vector ${\bf a} \in \mathbb{R}^m$, $m \in \mathbb{N} $. The {\bfi  inputs}  $ \left\{ {\bf z}_t\right\}_{t\in \Bbb Z}$ of the system, with ${\bf z}_t \in D_d$, will be in most cases infinite paths of a discrete-time stochastic process. We note that, unlike what we do in this paper, the term {\bfi  recurrent neural network}  is used sometimes in the literature to refer exclusively to state-space systems where the state map $F$ in \eqref{rc state eq} is neural network-like, that is, it is the composition of a nonlinear activation function with an affine function of the states and the input.

We shall focus on state-space systems of the type \eqref{rc state eq}-\eqref{rc readout eq} that determine an {\bfi  input/output} system. This happens in the presence of the so-called {\bfi  echo state property (ESP)}, that is, when for any  ${\bf z} \in (D_d)^{\mathbb{Z}}$ there exists a unique $\mathbf{y} \in (\mathbb{R} ^m)^{\mathbb{Z}}$ such that \eqref{rc state eq}-\eqref{rc readout eq} hold. In that case, we  talk about the {\bfi  state-space filter} $U ^F_h: (D_d)^{\mathbb{Z}}\longrightarrow  ({\mathbb{R}}^m)^{\mathbb{Z}} $ associated to the state-space system \eqref{rc state eq}-\eqref{rc readout eq}  defined by: 
\begin{equation*}
U  ^F _h({\bf z}):= \mathbf{y},
\end{equation*}
where ${\bf z} \in (D_d)^{\mathbb{Z}} $  and $\mathbf{y} \in (\mathbb{R}^m)^{\mathbb{Z}} $ are linked by \eqref{rc state eq} via the ESP. If the ESP holds at the level of the state equation \eqref{rc state eq}, we can define a {\bfi  state filter}  $U ^F: (D_d)^{\mathbb{Z}}\longrightarrow  (D_N)^{\mathbb{Z}} $ and, in that case, we have that
\begin{equation*}
U ^F_h:= h \circ U ^F. 
\end{equation*}
It is easy to show that state and state-space filters are automatically causal and time-invariant (see \cite[Proposition 2.1]{RC7}) and hence it suffices to work with their restriction $U ^F_h: (D_d)^{\mathbb{Z}_-}\longrightarrow  ({\mathbb{R}}^m)^{\mathbb{Z}_-} $ to semi-infinite inputs and outputs. Moreover, $U ^F_h $  determines a state-space {\bfi   functional} $H ^F_h: (D_d)^{\mathbb{Z}}\longrightarrow  {\mathbb{R}}^m$ as $H ^F_h({\bf z}):=U ^F_h({\bf z})_0 $, for all ${\bf z} \in (D _d)^{\mathbb{Z}_{-}} $ (the same applies to $U ^F $ and $H ^F $ when the ESP holds at the level of the state equation). In the sequel we use the symbol $\mathbb{Z}_{-}  $ to denote the negative integers including zero and $\mathbb{Z}^{-}  $ without zero.

The echo state property has received much attention in the context of the so-called {\bfi  echo state networks (ESNs)} \cite{Matthews:thesis, Matthews1993, Matthews1994, Jaeger04} (see, for instance,~\cite{jaeger2001, Buehner:ESN, Yildiz2012, zhang:echo, Wainrib2016, Manjunath:Jaeger, gallicchio:esp}). Sufficient conditions for the ESP to hold in general systems have been formulated in \cite{RC7, RC9, RC10}, in most cases assuming that the state map is a contraction in the state variable. We now recall a result (see \cite[Theorem 12]{RC9} and \cite[Proposition 1]{RC10}) that ensures the ESP as well as a continuity property of { the state filter which are important} for the sequel.
{ All along this paper and whenever in the presence of  Cartesian products (finite or infinite) of topological spaces, continuity will be considered with respect to the product topology, that is, the coarsest topology that makes continuous all the canonical projections onto the individual factors (see \cite[Chapter 2]{Munkres:topology} for details).} 
\begin{proposition}
\label{ESP for reservoir maps with compact target}
Let $F:  D _N \times  D _d \longrightarrow  D_N$  be a continuous state map such that  $D_N $ is a compact subset of ${\Bbb R}^N  $ and $F$ is a contraction on the first entry with constant $0<c<1 $, that is,
\begin{equation*}
\label{eq:Fcontractive} \|F({\bf x}_1,{\bf z})- F({\bf x}_2,{\bf z}) \| \leq c \|{\bf x}_1-{\bf x}_2\|, 
\end{equation*}
for all ${\bf x}_1,{\bf x}_2 \in D _N$, ${\bf z} \in D_d$.
Then, the  associated system { has} the echo state property for any input in $(D _d)^{\mathbb{Z}_{-}} $. The associated  filter $U ^F: (D _d)^{\mathbb{Z}_{-}}  \longrightarrow (D _N)^{\mathbb{Z}_{-}}$ is continuous with respect to the product topologies in $(D _d)^{\mathbb{Z}_{-}} $ and $(D _N)^{\mathbb{Z}_{-}} $. 
\end{proposition}
{ As we show in the next few paragraphs, a given filter admits non-unique representations which can be generated with maps between state spaces  resulting in  systems in general with  different memory and forecasting capacities. In the next proposition we show the important implication of the echo state property for our analysis. More specifically, we show that that whenever the map between state spaces is a system morphism, the target system of  the system morphism has the ESP and for the original system the existence of at least one solution for each input is guaranteed, then this solution is also unique or, equivalently,  this system also has ESP and, moreover, the filters associated to these two systems are identical. This result can be made even stronger for isomorphisms which  we will be using in our derivations.}

\paragraph{State-space morphisms.} {The state-space representations of a given filter, when they exist, are not necessarily unique. These different realizations can be generated using maps between  state spaces  that satisfy certain natural functorial properties that make them into morphisms in the category of state-space systems. Additionally, as we  see later on in Proposition \ref{properties of morphisms solutions}, morphisms encode information about the solution properties and the echo state property of the systems that are linked by them}. Consider the state-space systems determined by the two pairs $(F _i, h _i)$, $i \in \left\{1,2\right\}$, with $F _i: D_{N _i}\times D_ d\longrightarrow D_{N _i} $ and $h _i:D_{N _i} \longrightarrow \mathbb{R}^m $.
\begin{definition}
A map $f: D_{N_1} \longrightarrow D_{N_2}$ is a {\bfi  morphism} between the systems $(F_1, h_1)$ and $(F_2, h_2)$ whenever it satisfies the following two properties:
\begin{description}
    \item[(i)] {\bfi System equivariance:} $f(F_1({\bf x}_1, {\bf z})) = F_2(f({\bf x}_1), {\bf z})$, for all ${\bf x}_1 \in D_{N_1},$ and ${\bf z} \in D_d$.
    \item[(ii)] {\bfi Readout invariance:} $h_1({\bf x}_1) = h_2(f({\bf x}_1))$, for all ${\bf x}_1 \in D_{N_1}$.
\end{description}
\end{definition}

When the map $f$ has an inverse $f ^{-1} $  and this inverse is also a morphism between the systems determined by the pairs $(F_1, h_1)$ and $(F_2, h_2)$ we say that $f$ is a {\bfi system isomorphism} and that the systems $(F_1, h_1)$ and $(F_2, h_2)$ are {\bfi isomorphic}. Given a system $F_1: D_{N_1} \times D_d \longrightarrow D_{N_1}, h_1: D_{N_1} \longrightarrow \mathbb{R}^m$ and a bijection $f:D_{N_1} \longrightarrow D_{N_2}$, the map $f$ is a system isomorphism with respect to the system $F_2: D_{N_2} \times D_d \longrightarrow D_{N_2}, h_2: D_{N_2} \longrightarrow \mathbb{R}^m$ defined by 
\begin{align}
    F_2({\bf x}_2, {\bf z}) &:= f(F_1(f^{-1}({\bf x}_2), {\bf z})), \quad \text{for all} \quad {\bf x}_2 \in D_{N_2}, {\bf z} \in D_d, \label{isomorphic state map}\\
    h_2({\bf x}_2) &:= h_1(f^{-1}({\bf x}_2)), \quad \text{for all} \quad {\bf x}_2 \in D_{N_2}.\label{isomorphic readout map}
\end{align}

\begin{proposition}
\label{properties of morphisms solutions}
Let  $(F _i, h _i)$, $i \in \left\{1,2\right\}$, be two systems with $F _i: D_{N _i}\times D_ d\longrightarrow D_{N _i} $ and $h _i:D_{N _i} \longrightarrow \mathbb{R}^m $.  Let $f: D_{N_1} \longrightarrow D_{N_2}$ be a map. Then:
\begin{description}
\item[(i)] If $f$ is system equivariant and ${\bf x}^1 \in (D_{N_1})^{\mathbb{Z}_-}$ is a solution for the state system associated to $F_1$ and the input ${\bf z} \in (D_d)^{\mathbb{Z}_-}$, then so is $(f({\bf x}_t^1))_{t\in \mathbb{Z}_-} \in (D_{N_2})^{\mathbb{Z}_-}$ for the system associated to $F_2$ and the same input.
\item[(ii)] Suppose that the system determined by $(F_2, h_2)$ has the  echo state property and assume that the state system determined by $F_1$ has at least one solution for each element ${\bf z} \in (D_d)^{\mathbb{Z}_-}$. If $f$ is a morphism between $(F_1, h_1)$ and $(F_2, h_2)$, then $(F_1, h_1)$ has the echo state property and, moreover,
\begin{equation}\label{eq: morphism in proposition}
    U_{h_1}^{F_1} = U_{h_2}^{F_2}.
\end{equation}
\item [(iii)] If $f$ is a system isomorphism then the implications in the previous two points are reversible, that is, the indices $1$ and $2$ can be exchanged.
\end{description}
\end{proposition}

\paragraph{Input and output stochastic processes.} We now fix a probability space $(\Omega,\mathcal{A},\mathbb{P})$ 
on which all random variables are defined. The triple consists of the sample space $\Omega$, which is the set of possible outcomes, the $\sigma$-algebra $\mathcal{A}$ (a set of subsets of $\Omega$ (events)),  and a probability measure $\mathbb{P}:\mathcal{A}\longrightarrow [0,1]$. {The input  signal is modeled as a discrete-time stochastic process ${\bf Z} = ({\bf Z}_t)_{t \in \mathbb{Z}_-}$  taking values in $D_d \subset \mathbb{R}^d$. Moreover, we write ${\bf Z}(\omega)=({\bf Z}_t( \omega))_{t \in \mathbb{Z}_-}$ for each outcome $\omega \in \Omega$ to denote the realizations or sample paths of ${\bf Z}$ . Since  ${\bf Z}$ can be seen as a random sequence in $D_d \subset \mathbb{R}^d$, we write interchangeably ${\bf Z}:{\mathbb{Z}}_{-} \times \Omega \longrightarrow D_d$ and ${\bf Z}: \Omega\longrightarrow (D_d)^{{\mathbb{Z}}_{-}}$. The latter is by assumption measurable
with respect to the Borel $\sigma $-algebra induced by the product topology in $(D_d)^{{\mathbb{Z}}_{-}}$. In this paper we consider only memory reconstruction and forecasting information processing tasks and hence the target process which is commonly denoted in the literature by ${\bf Y}$, will always be a time backward or forward shifted version of the input time series process  ${\bf Z}$.}

We will most of the time work under {\bfi  stationarity} hypotheses. We recall that the discrete-time process ${\bf Z}: \Omega\longrightarrow (D_d)^{{\mathbb{Z}}_{-}}$ is stationary whenever $T _{-\tau}({\bf Z}) \stackon{=}{\mbox{{\tiny d}}} {\bf Z}$, for any $\tau\in \mathbb{Z}_{-}  $, where the symbol $\stackon{=}{\mbox{{\tiny d}}} $ stands for the equality in distribution and $T_{-\tau}: ({\Bbb R}^d)^{\mathbb{Z}_{-}} \longrightarrow ({\Bbb R}^d)^{\mathbb{Z}_{-}} $ is the {\bfi  time delay} operator defined by $T_{-\tau}({\bf z})_t:= {\bf z}_{t+ \tau} $ for any $t \in \mathbb{Z}_{-} $.

The definition of stationarity that we just formulated is usually known in the time series literature  as {\bfi  strict stationarity} \cite{BrocDavisYellowBook}. When a process ${\bf Z} $  has second-order moments, that is, ${\bf Z} _t \in L ^2(\Omega, {\Bbb R}^d)$, $t \in \mathbb{Z}_{-} $, then strict stationarity implies the so-called {\bfi  second-order stationarity}. We recall that  a square-integrable process  ${\bf Z}: \Omega\longrightarrow (D_d)^{{\mathbb{Z}}_{-}}$ is second-order stationary whenever {\bf (i)} there exists a constant ${ \boldsymbol{\mu}_{ Z}} \in {\Bbb R}^d $ such that ${\rm E}\left[{\bf Z} _t\right] ={ \boldsymbol{\mu}_{ Z}}$, for all $t \in \mathbb{Z}_{-}  $ ({\bfi  mean stationarity}) and {\bf (ii)} the autocovariance matrices $ \text{Cov} \left({\bf Z} _t , {\bf Z}_{t+h}\right)$  depend only on $h \in \Bbb Z  $ and not on $t \in \mathbb{Z}_{-} $ ({\bfi  autocovariance stationary}) and we can hence define the {\bfi   autocovariance function} {$\gamma: \Bbb Z \longrightarrow \mathbb{S}_d$ (with $\mathbb{S}_d$ the cone of positive semi-definite symmetric matrices of dimension $d$)} as  $\gamma(h):=\text{Cov} \left({\bf Z} _t , {\bf Z}_{t+h}\right)$, with $t\in \mathbb{Z}_{-}$ arbitrarily chosen so that $t+h \in \mathbb{Z}_{-} $. The autocovariance function necessarily satisfies $\gamma(h)= \gamma(-h)^{\top} $  \cite{BrocDavisYellowBook}. If $\mathbf{Z}  $ is mean stationary and condition {\bf (ii)} only holds for $h=0 $ we say that ${\bf Z}$ is {\bfi  covariance stationary}. Second-order stationarity and stationarity are only equivalent for Gaussian processes. If ${\bf Z} $ is autocovariance stationary and $\gamma(h)=0 $ for any non-zero $h \in \Bbb Z  $ then we say that ${\bf Z}$ is a {\bfi  white noise}.

\begin{corollary}
\label{stationary implies stationary}
Let $F:  D _N \times  D _d \longrightarrow  D_N$  be a  state map that satisfies the hypotheses of Proposition \ref{ESP for reservoir maps with compact target} or that, more generally, has the echo state property and the associated filter $U ^F: (D _d)^{\mathbb{Z}_{-}}  \longrightarrow (D _N)^{\mathbb{Z}_{-}}$ is continuous with respect to the product topologies in $(D _d)^{\mathbb{Z}_{-}} $ and $(D _N)^{\mathbb{Z}_{-}} $. If the input process ${\bf Z}: \Omega\longrightarrow (D_d)^{{\mathbb{Z}}_{-}}$ is stationary, then so is the {state} ${\bf X}:=U ^F({\bf Z}): \Omega \longrightarrow (D_N)^{\mathbb{Z}_{-}}$  as well as the joint processes $(T _{-\tau}({\bf X}), {\bf Z})$ and $({\bf X}, T _{-\tau}({\bf Z}))$, for any $\tau \in \mathbb{Z}_{-}  $.
\end{corollary}

{ In Proposition~\ref{properties of morphisms solutions} we showed how to design alternative state-representations of a given filter by using state morphisms. This freedom can be put at work by choosing representations that have specific technical advantages that are needed in a given situation. An important implementation example of this strategy is the next Proposition, where we show that if we have a state system and an input for which the output process is covariance stationary and the corresponding covariance matrix is non-singular, then an isomorphic system representation exists whose corresponding state process is standardized, that is,  the states have mean zero and covariance matrix equal to the identity. This standardization  leads to systems that are easier to handle in terms of the computation of memory and forecasting capacities, which is profusely exploited later on in the main results of the paper.}

\begin{proposition}[Standardization of state-space realizations]
\label{std rc io}
Consider a state-space system as in \eqref{rc state eq}-\eqref{rc readout eq} and suppose that the input process ${\bf Z}: \Omega\longrightarrow (D_d)^{{\mathbb{Z}}_{-}}$ is such that the associated state process ${\bf X}: \Omega\longrightarrow (D_N)^{{\mathbb{Z}}_{-}}$ is covariance stationary. Let $\boldsymbol{\mu}:= {\rm E}\left[\mathbf{X}_t\right] $  and suppose that the covariance matrix  $\Gamma_{{\bf X}} := {\rm Cov}(\mathbf{X}_t, \mathbf{X}_t)$ is non-singular.
Then, the map $f: \mathbb{R}^N \longrightarrow \mathbb{R}^N$ given by $f(\mathbf{x}):=\Gamma_{{\bf X}}^{-1/2} ( \mathbf{x} - \boldsymbol{\mu}) $ is a system isomorphism between the system \eqref{rc state eq}-\eqref{rc readout eq} and the one with state map \begin{equation}
\label{isomorphic state eq}
\widetilde{F}( \mathbf{x}, {\bf z}): =\Gamma _{{\bf X}}^{-1/2} \left(  F\left(\Gamma_{{\bf X}}^{1/2}  \mathbf{x} + \boldsymbol{\mu}, {\bf z}\right) -\boldsymbol{\mu} \right)
\end{equation}
and readout
\begin{equation}
\label{isomorphic readout eq}
\widetilde{h}( \mathbf{x}): = (\Gamma _{{\bf X}}^{-1/2} \mathbf{W})^\top \mathbf{x}+ {\bf W}^{\top} \boldsymbol{\mu} + \mathbf{a}.
\end{equation}
Moreover, the state process $\widetilde{{\bf X} } $ associated to the system $\widetilde{F} $  and the input ${\bf Z} $ is covariance stationary and
\begin{equation}
\label{mean and covariance new states}
{\rm E}[\widetilde{\mathbf{X}}_t] = {\bf 0}, \quad \mbox{and} \quad {\rm Cov}( \widetilde{\mathbf{X}}_t, \widetilde{\mathbf{X}}_t) = \mathbb{I}_N.
\end{equation}
\end{proposition}

\section{Memory and forecasting capacity bounds for stationary inputs}
\label{Memory and forecasting capacity bounds for stationary inputs}

The following definition extends the notion of memory capacity introduced in~\cite{Jaeger:2002} to general nonlinear systems and to input signals that are stationary but not necessarily time-decorrelated.

\begin{definition}
\label{memory forecasting capacities}
Let  ${\bf Z}: \Omega\longrightarrow D^{{\mathbb{Z}}_{-}}$, $D \subset \mathbb{R}$, be a variance-stationary input and let $F$ be a state map that has the echo state property with respect to the paths of ${\bf Z} $. Assume, moreover, that the associated state process ${\bf X}: \Omega\longrightarrow (D_N)^{{\mathbb{Z}}_{-}}$ defined by ${\bf X}_t:=U ^F({\bf Z})_t  $ is covariance stationary, as well as the joint processes $(T _{-\tau}({\bf X}), {\bf Z})$ and $({\bf X}, T _{-\tau}({\bf Z}))$, for any $\tau \in \mathbb{Z}_{-}  $. We define the $\tau$-lag {\bfi  memory capacity} ${\rm MC} _\tau $ (respectively, {\bfi  forecasting capacity} ${\rm FC} _\tau $) of $F$ with respect to ${\bf Z} $ as:
\begin{eqnarray}
{\rm MC} _\tau &:=& 1- \frac{1}{\text{\rm Var} \left(Z _t\right)}  \stackunder{{\rm min}}{\stackunder{\scriptstyle {\bf W} \in \mathbb{R}^N }{\scriptstyle a \in \mathbb{R} }} {\rm E}\left[\left(\left(T _{-\tau} {\bf Z}\right)_t- {\bf W} ^{\top} U ^F({\bf Z})_t-a\right)^2\right],\label{definition memory t}\\
{\rm FC} _\tau &:=& 1- \frac{1}{\text{\rm Var} \left(Z _t\right)}  \stackunder{{\rm min}}{\stackunder{\scriptstyle {\bf W} \in \mathbb{R}^N }{\scriptstyle a \in \mathbb{R} }} {\rm E}\left[\left( {Z}_t- {\bf W} ^{\top} U ^F(T _{-\tau}({\bf Z}))_t-a\right)^2\right].\label{definition forecasting t}
\end{eqnarray}
The {\bfi total memory capacity} ${\rm MC}$ (respectively, {\bfi  total forecasting capacity} ${\rm FC}$) of $F$ with respect to ${\bf Z} $ is defined as:
\begin{equation}
\label{for and mem capacities}
{\rm MC} := \sum _{\tau \in \mathbb{Z}_{-}} {\rm MC} _\tau, \qquad
{\rm FC} := \sum _{\tau \in \mathbb{Z}^{-}} {\rm FC} _\tau.
\end{equation}
\end{definition}

Note that, by Corollary \ref{stationary implies stationary}, the conditions of this definition are met when, for instance,
$U ^F$ is continuous with respect to the product topologies and the input process ${\bf Z}$ is stationary. 

{ The optimization problems appearing in the definitions \eqref{definition memory t} and \eqref{definition forecasting t} of the memory and forecasting capacities can be explicitly solved when the state covariance matrix $\Gamma_{{\bf X}} := {\rm Cov}(\mathbf{X}_t, \mathbf{X}_t)$ is invertible. We refer to this situation as the {\it regular case}. These solutions are provided in the following lemma.}
\begin{lemma}
\label{rewrite with covs}
In the conditions of Definition \ref{memory forecasting capacities} and if the covariance matrix  $\Gamma_{{\bf X}} := {\rm Cov}(\mathbf{X}_t, \mathbf{X}_t)$ is invertible then, for any $\tau \in \mathbb{Z}_{-} $:
\begin{equation}
\label{capacity formulas when invertible}
{\rm MC} _\tau = \frac{{\rm Cov} \left(Z_{t+ \tau}, {\bf X} _t\right) \Gamma_{{\bf X}} ^{-1} {\rm Cov} \left({\bf X} _t, Z_{t+ \tau}\right)}{\text{\rm Var} \left(Z _t\right)}, \qquad
{\rm FC} _\tau =\frac{{\rm Cov} \left(Z_{t}, {\bf X} _{t+ \tau}\right) \Gamma_{{\bf X}} ^{-1} {\rm Cov} \left({\bf X} _{t+ \tau}, Z_{t}\right)}{\text{\rm Var} \left(Z _t\right)}.
\end{equation}
\end{lemma}
{The availability of the closed-form solutions in the expressions \eqref{definition memory t} and \eqref{definition forecasting t} make the computation of capacities much easier. This will become particularly evident in the next section devoted to linear systems.  We emphasize that even for those simpler system specifications, the non-invertibility of the associated covariance matrices of states leads to technical difficulties.  The same holds even to a greater extent for nonlinear systems. Some of those problems can be handled by using equivalent state-space representations. The next result shows that new representations obtained out of linear injective system morphisms leave invariant the capacities and hence can be used to produce systems with more technically tractable properties.} This result will be used later on in Section \ref{The linear case} when we study the memory and forecasting capacities {of linear systems in the singular case}. 
\begin{lemma}
\label{capacity invariant linear morphism}
Let ${\bf Z}  $ be a variance-stationary input and let $F_2: D_{N _2}\times D \longrightarrow D_{N _2} $ be a state map that satisfies the conditions of Definition \ref{memory forecasting capacities}. Let $F_1: D_{N _1}\times D \longrightarrow D_{N _1} $ be another state map that has at least one solution for each ${\bf z} \in D ^{\mathbb{Z}_{-}} $ and let $f: \mathbb{R}^{N_1}\longrightarrow\mathbb{R}^{N_2}$ be an injective linear system equivariant map between $F _1  $  and $F _2$. Then, the memory and forecasting capacities of $F _1 $ with respect to ${\bf Z} $ are well-defined and coincide with those of $F _2 $ with respect to ${\bf Z} $.
\end{lemma}

{The next theorem is the first main contribution} of this paper and generalizes the bounds formulated in~\cite{Jaeger:2002} for the total memory capacity of an echo state network in the presence of independent inputs to general state systems with second-order stationary inputs {and invertible state covariance matrices. We show that both the total memory and  forecasting capacities are nonnegative and that upper bounds can be formulated that are fully determined by the behavior of the autocovariance or the spectral density functions of the input and the dimensionality of the state space.}

\begin{theorem}
\label{capacity bounds}
Suppose that we are in the conditions of Definition \ref{memory forecasting capacities} and that the covariance matrix  $\Gamma_{{\bf X}} := {\rm Cov}(\mathbf{X}_t, \mathbf{X}_t)$ is non-singular. 
\begin{description}
\item [(i)] For any $\tau \in \mathbb{Z}_{-} $:
\begin{equation}
\label{bounds fct and mct}
0 \leq {\rm MC} _\tau\leq 1 \quad \mbox{and} \quad 0 \leq {\rm FC} _\tau\leq 1.
\end{equation}
\item [(ii)] Suppose that, additionally, the input process $ {\bf Z}  $ is second-order stationary with  autocovariance function $\gamma: \Bbb Z \longrightarrow \mathbb{R} $, and that for any $L\in \mathbb{Z}_{-}$ the symmetric matrices $H ^L \in \mathbb{M}_{-L+1} $  determined by $H ^L _{ij}:= \gamma (|i-j|) $ are invertible. Then, if we use the symbol ${\rm C}$ to denote both ${\rm MC} $ and  ${\rm  FC}$ in \eqref{for and mem capacities}, we have:
\begin{equation}
\label{bounds memory cap}
0\leq  {\rm C}\leq \frac{N}{\gamma (0)}\rho(H)\leq N \left(1+ \frac{2}{\gamma (0)}\sum _{j=1}^{\infty}|\gamma ( j) |\right),
\end{equation}
where $\rho(H):=\lim\limits_{L \rightarrow - \infty} \rho(H^L)$, with $\rho(H^L)$ the spectral radius of $H ^L $.
\item [(iii)] In the same conditions as in part {\bf (ii)}, suppose that, additionally, the autocovariance function $\gamma$ is absolutely summable, that is, $\sum _{j=- \infty}^{\infty}|\gamma (j)|< +\infty $. In that case, the spectral density $f :[- \pi, \pi]\rightarrow  \mathbb{R}  $  of ${\bf Z}  $ is well-defined and given by
\begin{equation}
\label{spectral density formula}
f (\lambda):= \frac{1}{2 \pi}\sum_{ n \in \Bbb Z} e^{-in \lambda}\gamma (n),
\end{equation}
and we have that
\begin{equation}
\label{inequality summable version}
0\leq {\rm C}\leq \frac{2 \pi N}{\gamma (0)}M _f\leq N \left(1+ \frac{2}{\gamma (0)}\sum _{j=1}^{\infty}|\gamma (j)|\right),
\end{equation}
where $M _f:= \max_{\lambda \in [- \pi, \pi]} \left\{f (\lambda)\right\}. $
\end{description}
\end{theorem}
{ We emphasize that the results in this theorem hold for a general class of nonlinear recurrent neural networks with linear readouts both with trainable or with  randomly generated neuron weights (reservoir computing). Additionally, it can be used as a  tool in the design of the network architecture when the autocovariance structure of the input is known. These bounds show in passing that memory and forecasting capacities are determined not only by a system but also to a great extent by the memory of the input process itself. }

\medskip

This {general} result and its proof are used in the next corollary to recover  the {total} memory capacity bounds proposed in \cite{Jaeger:2002} when using independent inputs and to show that, in that case, the {total} forecasting capacity is {always} zero.

\begin{corollary}[\cite{Jaeger:2002}]
\label{jaeger bound}
In the conditions of Theorem \ref{capacity bounds}, if the inputs $\left\{Z _t\right\}_{t \in \mathbb{Z}_{-}}$ are independent, then
\begin{equation}
\label{bounds in independent case}
0\leq {\rm MC}\leq N \quad \mbox{and} \quad {\rm FC}=0.
\end{equation}
\end{corollary}

{The proofs of the previous two results, which can be found in the appendices, shed some light on the relative values of the forecasting and memory capacities, as well as on the quality of the common bounds in \eqref{bounds memory cap} and \eqref{inequality summable version}. Indeed, these estimates are obtained by finding upper bounds for the norms of the orthogonal projections (in the $L ^2 $ sense) of the state at a given time on the vector space generated by all the inputs fed into the system up until that point in time, in the case of the memory capacity and, for the forecasting capacity, by the inputs that will be fed in the future. The built-in causality of state-space filters implies that the state has a functional dependence exclusively on past inputs and hence its projection onto future inputs becomes non-trivial only via dependence phenomena in the input signal. This fact carries in its wake that, typically, even in the presence of strongly autocorrelated input signals, the projection of the state vector onto past inputs produces larger vectors (in norm) than onto future ones. The bounding mechanism used in the proof (see \eqref{bounding by projection}) is not able to take this fact into account, as that would entail using specific knowledge on the functional form of the filter, which is something that we avoided in the pursuit of generic bounds that are common to all state-space systems with a given dimension. The price to pay for this degree of generality is that the bounds will be closer to the memory than to the forecasting capacities and hence will be sharper for the former than for the latter. Later on in Section \ref{Numerical illustration} we illustrate these facts with a numerical example and we additionally discuss the sharpness question, or rather the lack of it, in the presence of dependent inputs.}

\section{The memory and forecasting capacities of linear systems}
\label{The linear case}

When the state equation \eqref{rc state eq} is linear and has the echo state property, both the memory and forecasting capacities in \eqref{for and mem capacities} can be explicitly written down in terms of the equation parameters provided that the invertibility hypothesis on the covariance matrix of the states holds. This case has been studied for independent inputs and randomly generated linear systems in \cite{Jaeger:2002, linearESN} and, more recently, in \cite{marzen:capacity} for more general correlated inputs and diagonalizable linear systems. It is in this paper that it has been pointed out for the first time how different linear systems that maximize forecasting and memory capacities may be. 

{This section contains the second main contribution of the paper}. We split it in two parts. In the first one we handle what we call the regular case in which we assume the invertibility of the covariance matrix of the states. In the second one we see how, using {system morphisms} of the type introduced in Lemma \ref{capacity invariant linear morphism}, we can reduce the general singular case to the regular one. This approach allows us to prove that {\it when the inputs are independent, the memory capacity of a linear system with independent inputs coincides with the rank of its associated controllability or reachability matrix}. { A rigorous proof of this result is, to our knowledge, not available in the literature}. This statement is a generalization of the fact, already established in \cite{Jaeger:2002}, that when the rank of the controllability matrix is maximal, then the linear system has maximal capacity, that is, its capacity coincides with the dimensionality of its state space. Different configurations that maximize the rank of the controllability matrix have been recently studied in \cite{Aceituno2017, Verzelli2020}. { We find the results in this section useful also from the applications point of view as they allow the design of linear recurrent networks with an exact pre-specified memory capacity.}

\subsection{The regular case}
\label{The regular case}

The {explicit} capacity formulas that we state in the next result only require the stationarity of the input and the invertibility of the covariance matrix of the { states}. In this case we consider fully infinite inputs, that is, we  work with a stationary input process ${\bf Z}: \Omega \longrightarrow \mathbb{R}^{\mathbb{Z}} $ that has second-order moments and  an autocovariance function $\gamma: \Bbb Z \longrightarrow \mathbb{R} $ for which we  construct the semi-infinite (respectively, doubly infinite) symmetric positive semi-definite Toeplitz matrix $H _{ij}:=  \gamma (|i-j|) $, $i,j \in \mathbb{N} ^+ $, (respectively, $\overline{H} _{ij}:=  \gamma (|i-j|) $, $i,j \in \Bbb Z$).

\begin{proposition}
\label{the linear case with bs}
Consider the linear state system determined by the linear state map $F: {\Bbb R}^N \times \mathbb{R} \longrightarrow{\Bbb R}^N $ given by
\begin{equation}
\label{definition linear system}
F(\mathbf{x}, z):= A \mathbf{x}+ \mathbf{C}z, \quad \mbox{where} \quad \mathbf{C} \in {\Bbb R}^N,  A \in \mathbb{M}_N,\,\mbox{and}\ \, \vertiii{A}=\sigma_{{\rm max}}(A)<1,
\end{equation}
{with $\sigma_{{\rm max}}(A)$ the largest singular value of the matrix $A$  (usually referred to as {\bfi  connectivity matrix}).}
Let $D\subset \mathbb{R}$ be compact and consider a zero-mean stationary input process ${\bf Z}: \Omega \longrightarrow D^{\mathbb{Z}} $ that has second-order moments and an absolutely summable autocovariance function $\gamma: \Bbb Z \longrightarrow \mathbb{R} $. Suppose also that the associated spectral density $f$ satisifies that $f (\lambda)\geq 0 $, $\lambda \in [- \pi, \pi]$, and that $f (\lambda)=0 $ holds only in at most a countable number of points. 

Then $F$ has the echo state property and the associated filter $U^{A, \mathbf{C}} $ is such that its output ${\bf X}:=U ^{A, \mathbf{C}}({\bf Z}): \Omega \longrightarrow (D_N)^{\mathbb{Z}_{-}}$  as well as the joint processes $(T _{-\tau}({\bf X}), {\bf Z})$ and $({\bf X}, T _{-\tau}({\bf Z}))$, for any $\tau \in \mathbb{Z}_{-}  $, are stationary and ${\bf X}  $ is covariance stationary.  Suppose that the covariance matrix $\Gamma_{{\bf X}} := {\rm Cov}(\mathbf{X}_t, \mathbf{X}_t)$  is non-singular. Then, 
\begin{description}
\item [(i)] Consider the $N$ vectors $\mathbf{B} _1, \ldots, \mathbf{B} _N   \in \ell^2_+(\mathbb{R}) $ defined by 
\begin{equation}
\label{definition b1}
{B}_i^j :=\left(\Gamma_{{\bf X}}^{-1/2}\sum _{k=0}^{\infty}A ^k \mathbf{C}H_{k+1,j}^{1/2}\right)_i, \quad j \in \mathbb{N}^+, \enspace { i =1,\ldots, N},
\end{equation}
where the square root matrices are computed via orthogonal diagonalization. These entries are all finite and form vectors that constitute an orthonormal set in $\ell^2_+(\mathbb{R}) $. The {total} memory capacity {\rm MC} can be written as
\begin{equation}
\label{expression memory linear case}
{\rm MC}= \frac{1}{\gamma(0)}\sum_{i=1}^N \langle \mathbf{B} _i, H \mathbf{B}_i \rangle _{\ell^2}.
\end{equation}
\item [(ii)] Consider the $N$ vectors $ \mathbf{B} _1, \ldots, \mathbf{B} _N   \in \ell^2(\mathbb{R}) $ defined by 
\begin{equation}
\label{definition b2}
{B}_i^j :=\left(\Gamma_{{\bf X}}^{-1/2}\sum _{k=0}^{\infty}A ^k \mathbf{C}\overline{H}_{-k,j}^{1/2}\right)_i, \quad j \in \Bbb Z, \enspace { i =1,\ldots, N}.
\end{equation}
These vectors form an orthonormal set in $\ell^2(\mathbb{R})$ and the forecasting capacity {\rm FC} can be written as
\begin{equation}
\label{expression forecasting linear case}
{\rm FC}= \frac{1}{\gamma(0)}\sum_{i=1}^N 
\left\|\mathbb{P}_{\mathbb{Z}^+}\left(\overline{H}^{1/2}\mathbf{B} _i\right)\right\|_{\ell^2}^2,
\end{equation}
where $\mathbb{P}_{\mathbb{Z}^+}: \ell^2(\mathbb{R})\longrightarrow \ell^2(\mathbb{R}) $ is the projection that sets to zero all the entries with non-positive index.
\end{description}
\end{proposition}

When the inputs are second-order stationary and not autocorrelated ($ {\bf Z} $ is a white noise) then the formulas in the previous result can be used to give, for the linear case, a more informative version of Corollary \ref{jaeger bound}, without the need to invoke input independence.

\begin{corollary}
\label{jaeger bound for linear}
Suppose that we are in the hypotheses of Proposition \ref{the linear case with bs} and that, additionally, the input process ${\bf Z}: \Omega \longrightarrow D^{\mathbb{Z}} $ is a white noise, that is, the autocovariance function $\gamma$ satisfies that $\gamma(h)=0 $, for any non-zero $h\in \Bbb Z$. Then:
\begin{equation}
\label{traditional capacities linear}
{\rm MC}=N \quad \mbox{and } \quad {\rm FC}=0.
\end{equation}
\end{corollary}

An important conclusion of this corollary is that linear systems with white noise inputs that have a non-singular state covariance matrix $\Gamma_{{\bf X}}$ automatically have maximal memory capacity. This makes important the characterization of the invertibility of $\Gamma_{{\bf X}}$ in terms of the  parameters $A \in \mathbb{M}_N $ and $\mathbf{C} \in \mathbb{R}^N$ in \eqref{definition linear system}.  { The following proposition provides such a  characterization, which has serious practical implications at the time of designing linear recurrent networks, and establishes} a connection between the invertibility of $\Gamma_{{\bf X}}$ and Kalman's characterization of the controllability of a linear system \cite{Kalman2010}.
\begin{proposition}
\label{capacities and kalman}
Consider the linear state system introduced in \eqref{definition linear system}. Suppose that the input process ${\bf Z}: \Omega \longrightarrow D^{\mathbb{Z}} $ is a white noise and that the connectivity matrix $A \in \mathbb{M}_N $ is diagonalizable. Let $\sigma(A)= \left\{\lambda_1, \ldots, \lambda_N\right\} $ be the spectrum of $A$ and let $\{ \mathbf{v}_1, \ldots, \mathbf{v}_N\}  $ be an eigenvectors basis. 
\begin{description}
\item [(i)] The state covariance matrix $\Gamma_{{\bf X}}= {\rm Cov}(\mathbf{X}_t, \mathbf{X}_t)= \gamma(0)\sum_{j=0}^{\infty}A ^j \mathbf{C}\mathbf{C} ^{\top} \left(A ^j\right)^{\top}$ is non-singular if and only if all the eigenvalues in $\sigma(A) $ {are} distinct.
\item [(ii)] The vectors $\left\{A \mathbf{C}, A ^2\mathbf{C}, \ldots, A ^N \mathbf{C}\right\}$ form a basis of $ \mathbb{R}^N$ if and only if all the eigenvalues in $\sigma(A) $ are distinct and non-zero and in the linear decomposition $\mathbf{C}=\sum_{i =1}^N c _i \mathbf{v} _i $,  the coefficients $c _i $, $i \in \left\{1, \ldots, N\right\}$ are all non-zero.
\item [(iii)] The conditions in the previous point are equivalent to the {\bf Kalman controllability condition}: the vectors $\left\{\mathbf{C}, A \mathbf{C}, \ldots, A ^{N-1} \mathbf{C}\right\}$ form a basis of $ \mathbb{R}^N$ together with the condition that all the eigenvalues in $\sigma(A) $ are non-zero.
\end{description}
\end{proposition}
{ It has been shown in \cite{Jaeger:2002} that the controllability condition is equivalent to maximum memory capacity.  
 The approach followed in this proposition will allow us  in Theorem \ref{expression capacity linear system} to generalize this statement by proving that the memory capacity {\it equals} the rank of the controllability matrix (introduced in detail later on in the text).}

\subsection{The singular case}

In the following paragraphs we study the situation in which the covariance matrix of the states $\Gamma_{{\bf X}}$ with white noise inputs is not invertible. All the results formulated in the paper so far, in particular the capacity formulas in \eqref{capacity formulas when invertible}, are not valid anymore in this case. {However, it is well-known that the non-invertibility of the  covariance matrix of the states process for systems with high state space dimensionality is a frequent  issue. In the reservoir computing literature this problem is usually overcome at the time of training via the use of spectral regularization techniques, like for instance the Tikhonov regularized regressions. In this paper we adopt a different strategy to tackle this problem, namely we use the idea introduced in Lemma~\ref{capacity invariant linear morphism} of using system morphisms that leave capacities invariant.} More specifically, we show that whenever we are given a linear system whose covariance matrix $\Gamma_{{\bf X}}$ is not invertible, there exists another linear system defined in a dimensionally smaller state space that generates the same filter and hence has the same capacities but, unlike the original system, this smaller one has an invertible covariance matrix. This feature allows us to use for this system some of the results in previous sections and, in particular {\it  to compute its memory capacity in the presence of independent inputs} that, as we establish in the next theorem, {\it coincides with the rank of the controllability matrix}.

\begin{theorem}
\label{expression capacity linear system}
Consider the linear system $F(\mathbf{x},z):=A \mathbf{x}+ \mathbf{C}z $  introduced in \eqref{definition linear system} and suppose that the input process ${\bf Z}: \Omega \longrightarrow D^{\mathbb{Z}} $ is a strictly stationary white noise. {Let $R(A, \mathbf{C}):= \left(\mathbf{C}| A \mathbf{C}| \cdots| A ^{N-1} \mathbf{C}\right) $ be the controllability matrix of the linear system.}
\begin{description}
\item [(i)] {Then it holds that }
\begin{equation}
\label{kernel of gammax}
\ker \Gamma_{\mathbf{X}}=\ker R(A, \mathbf{C}) ^{\top}.
\end{equation}
\item [(ii)] Let $X:= {\rm span}\left\{\mathbf{C}, A \mathbf{C}, \ldots, A ^{N-1} \mathbf{C}\right\} $ and let $r:=\dim (X)={\rm rank} \,R(A, \mathbf{C})$. Let $V\subset \mathbb{R}^N$ be a vector subspace such that $\mathbb{R}^N=X\oplus V $ and let $i _X: X \hookrightarrow \mathbb{R}^N  $ and  $\pi_X: \mathbb{R}^N \longrightarrow X $ be the injection and the projection associated to this splitting, respectively. Then, the linear system $\overline{F}: X \times D \longrightarrow X $ defined by 
\begin{equation}
\label{equivalent linear system}
\overline{F}(\overline{\mathbf{x}},z):= \overline{A} \overline{\mathbf{x}}+ \overline{\mathbf{C}} z \  \mbox{and determined by $\overline{A}:= \pi _XA i_X $ and $ \overline{\mathbf{C}}= \pi_X(\mathbf{C}) $,}
\end{equation}
 is well-defined and has the echo state property.
\item [(iii)] { In the notation of part {\bf{(ii)}}, the map} $i _X: X \hookrightarrow \mathbb{R}^N  $ is an injective linear system equivariant map between $\overline{F}  $  and $F $. 
\item [(iv)]  { In the notation of part {\bf{(ii)}}}, let $\overline{{\bf X} }: \Omega \longrightarrow X ^{\mathbb{Z}} $ be the output of the filter determined by the state-system $\overline{F} $. Then,
\begin{equation}
\label{invertibility of the new system}
{\rm rank}\, R(\overline{A}, \overline{\mathbf{C}})={{\rm rank} \,R(A, \mathbf{C}),}
\end{equation}
and if $\overline{A} $ is diagonalizable with non-zero eigenvalues then $\Gamma_{\overline{{\bf X}}}:= {\rm Cov}(\overline{\mathbf{X}}_t, \overline{\mathbf{X}}_t) $ is invertible.
\item [(v)] If $\overline{A} $ is diagonalizable with non-zero eigenvalues then the memory ${\rm MC} $ and forecasting ${\rm FC}$ capacities of $F$  with respect to ${\bf Z} $  are given by
\begin{equation}
\label{traditional capacities linear general}
{\rm MC}= {\rm rank}\,R(A, \mathbf{C})= \dim \left({\rm span} \left\{\mathbf{C}, A \mathbf{C}, \ldots, A^{N-1}\mathbf{C}\right\} \right)\quad \mbox{and } \quad {\rm FC}=0.
\end{equation}
\end{description}
\end{theorem}
{ The statement in part {\bf (v)} provides a generalization of the statement  in \cite{Jaeger:2002} that establishes the equivalence between  controllability (${\rm rank}\,R(A, \mathbf{C})=N $) and maximum memory capacity (${\rm MC}=N$).  More specifically, our statement shows that  the memory capacity {\it equals} the rank of the controllability matrix of the linear system. This result has far reaching implications for the applications  of recurrent linear networks with either fully trainable or just randomly generated weights.  Given a precise computational task at hand, a learner can use the rank of the controllability matrix  in order to construct a network with a prescribed memory  capacity. In particular, in the reservoir computing community the result in \cite{Jaeger:2002} has deserved much attention. Since the the controllability condition was known to be equivalent to maximal memory capacity, many attempts have been made trying to propose a strategy to generate  random reservoirs which would have maximal expected controllability matrix rank \cite{Rodan2011, Aceituno2017, Tino2019, Verzelli2020}. Our results show that the same work can be done in those cases when one is interested in constructing random reservoirs with a required  controllability matrix rank and that, as we proved, amounts to  the memory capacity of the system.}

\subsection{Numerical illustration}
\label{Numerical illustration}

{An important consequence of part {\bf (v)} in Theorem \ref{expression capacity linear system} is that the capacity bounds in \eqref{bounds memory cap} and \eqref{inequality summable version} are sharp in the presence of independent inputs or, equivalently, that the bounds in Corollary \ref{jaeger bound} are sharp. Indeed, by \eqref{traditional capacities linear general}, any linear system with independent inputs whose controllability matrix has maximal rank has full memory capacity  equal to its dimension, and hence it achieves the upper bound in Corollary \ref{jaeger bound}. A natural question that arises is if this sharpness remains valid for correlated inputs. Even though we were not able to formulate general conditions that would ensure that fact, the following paragraphs contain a numerical illustration that give indications of what the situation may be. Indeed, we demonstrate that, in general, {\it the controllability condition does not ensure anymore the sharpness of the bounds \eqref{bounds memory cap} and \eqref{inequality summable version} with correlated inputs} neither for memory capacities nor for forecasting capacities. The latter is a consequence of the arguments in the paragraph after Corollary \ref{jaeger bound}.}

\begin{figure}[h]
\label{fig:capacities_figure}
\vspace{0cm}
\includegraphics[scale=0.32, angle=0]{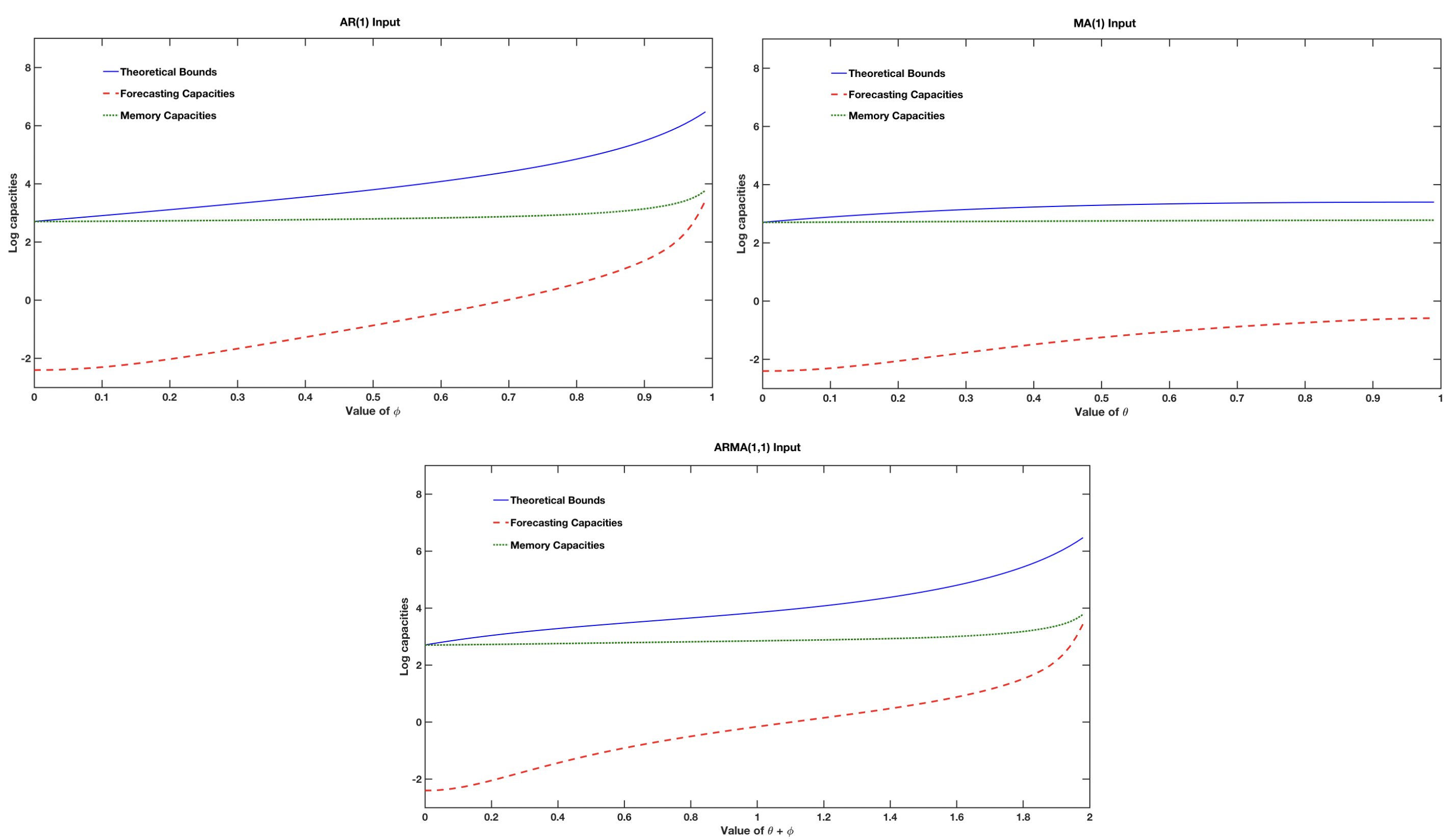}
\caption{Numerically computed memory and forecasting capacities of a linear system with full rank  controllability matrix and AR(1), MA(1), and ARMA(1,1) inputs. The curves depict the behavior of the numerically computed capacities and the bounds in \eqref{bounds memory cap} when the input model parameters are varied. These results show that, in this case, the theoretical bounds are only sharp for independent inputs.}
\end{figure}

{The panels in Figure \ref{fig:capacities_figure} show (in logarithmic scale) numerically computed memory and forecasting capacities, as well as the bounds in \eqref{bounds memory cap} based on the spectral radius $\rho(H)  $, for a linear system as in \eqref{definition linear system}. In this experiment we chose $N=15$ and a connectivity matrix $A$ (spectral radius equal to $0.9 $) and an input vector $\mathbf{C} $ such that  ${\rm rank}\,R(A, \mathbf{C})=N=15 $. The resulting system has hence memory capacity equal to $15  $ in the presence of independent inputs.} 

{This system has been then presented with three different types of autocorrelated inputs that are realizations of AR(1), MA(1), and ARMA(1,1) processes (see, for instance, \cite{BrocDavisYellowBook} for details on these models) driven by independent standard normal innovations. We denote (as in the figure) by $\phi $ and $\theta $ the autoregressive and the moving-average coefficients needed in the specification of these models. The top two panels in Figure \ref{fig:capacities_figure} have been obtained by varying the values of $\phi $ (for the AR(1) case) and of $\theta$ (for the MA(1) case) between $0$ and $1$. In the one at the bottom we took $\phi $ equal to $\theta  $ and we then varied them simultaneously between $0$ and $1$.} 

{The curves in the figures show how the bounds and the capacities evolve as a function of those parameters. The capacities have been computed using the definition in \eqref{for and mem capacities} and the formulas \eqref{capacity formulas when invertible} where we truncated the infinite sum at the value $\tau=250 $ and the covariances where empirically estimated using realizations of  length ten thousand.
In both cases, the values $\phi=0 $  and $\theta=0 $ correspond to the independent inputs case and the figures show how then the theoretical bounds correspond to the actual memory capacity of the system equal to $15 $. As soon as both parameters are non-zero we see in the figures a monotonous increase in the memory and forecasting capacities of the system, which is specially visible when it comes to the relation between the autoregressive parameter $\phi $ and the forecasting capacity. The figures also show that, as we anticipated, the theoretical bounds are strictly above the memory capacity even though we are in the presence of a system with full rank controllability matrix. This shows, in passing, that the results in Theorem \ref{expression capacity linear system} do not automatically extend to the dependent inputs case.
}

\section{Conclusions}
\label{Conclusions}

In this paper we have studied memory and forecasting capacities of generic nonlinear recurrent networks with respect to arbitrary stationary inputs that are not necessarily independent. In particular, we have stated {\it upper bounds for total memory and forecasting capacities in terms of the dimensionality of the network and the autocovariance of the inputs} that generalize those formulated in~\cite{Jaeger:2002}  for independent inputs.

The approach followed in the paper is particularly advantageous for linear networks for which explicit expressions can be formulated  for both  capacities. 
In the {\it classical} linear case with independent inputs, we have proved that {\it the memory capacity of a linear recurrent network with independent inputs is given by the rank of its controllability matrix}. This explicit and readily computable characterization of the memory capacity of those networks generalizes a well-known relation between maximal capacity and Kalman's controllability condition, formulated for the first time in  \cite{Jaeger:2002}. This is, to our knowledge, the first rigorous proof of the relation between network memory and the rank of the controllability matrix, that has been for a long time  part of the reservoir computing folklore.

The results in this paper suggest links between controllability and memory capacity for nonlinear recurrent systems that will be explored in forthcoming works.

\bigskip

\section{Appendices}
\label{Appendices}

\subsection{Proof of Proposition~\ref{properties of morphisms solutions}}
 \noindent\textbf{(i)} By hypothesis ${\bf x}_t^1 = F_1({\bf x}_{t-1}^1, {\bf z}_t)$, for all $ t \in \mathbb{Z}_{-}$. By system equivariance, 
$$
f({\bf x}_t^1) = f(F_1({\bf x}_{t-1}^1, {\bf z}_t))= F_2(f({\bf x}_{t-1}^1), {\bf z}_t), \quad \mbox{\,for all $ t \in \mathbb{Z}_{-}$.}
$$ 
as required.

\medskip

\noindent \textbf{(ii)} In order to show that \eqref{eq: morphism in proposition} holds, it suffices to prove that given any ${\bf z} \in (D_d)^{\mathbb{Z}_-}$, any solution $({\bf y}^1, {\bf x}^1) \in \left( \mathbb{R}^{m}\right)^{\mathbb{Z}_-} \times \left( D_{N _1}\right)^{\mathbb{Z}_-}$ of $(F_1, h_1)$ associated to ${\bf z}$, with ${\bf y}^1:= \left(h _1 ({\bf x}^1 _t)\right)_{t \in \mathbb{Z}_{-}}$, which exists by the hypothesis on $F _1 $,  coincides with the unique solution $U_{h_2}^{F_2}({\bf z})$ for the system $(F_2, h_2)$. Indeed, for any $t \in \mathbb{Z}_-$,
$$
{\bf y}_t^1 = h_1(F_1({\bf x}_{t-1}^1, {\bf z}_t)) = h_2(f(F_1({\bf x}_{t-1}^1, {\bf z}_t))) =h_2(F_2(f({\bf x}_{t-1}^1), {\bf z}_t)).$$
Here, the second equality follows from  readout invariance and the third one from the system equivariance. This implies that $\left({\bf y}^1, (f({\bf x}_t^1))_{t \in \mathbb{Z}_-}\right)$ is a solution of the system determined by $(F_2, h_2)$ for the input ${\bf z} $. By hypothesis, $(F_2, h_2)$ has the echo state property and hence ${\bf y}^1 = U_{h_2}^{F_2}({\bf z})$ and since ${\bf z} \in (D_d)^{\mathbb{Z}_-}$ is arbitrary, the result follows. Part {\bf (iii)} is straightforward.
\quad $\blacksquare$

\subsection{Proof of Corollary~\ref{stationary implies stationary}}
\noindent The continuity (and hence the measurability) hypothesis on $U^F $ proves that ${\bf X}:=U ^F({\bf Z}) $ is stationary (see \cite[page 157]{Kallenberg2002}). The joint processes $(T _{-\tau}({\bf X}), {\bf Z})$ and $({\bf X}, T _{-\tau}({\bf Z}))$ are  also stationary as they are the images of $ {\bf Z} $ by the measurable maps $(T _{-\tau} \circ U ^F) \times \mathbb{I}_{(\mathbb{R}^d)^{\mathbb{Z}_{-}}}$ and $U ^F \times T _{-\tau}$, respectively. \quad $\blacksquare$

\subsection{Proof of Proposition~\ref{std rc io}}
\noindent Since by hypothesis the matrix $\Gamma _{{\bf X}} $  is invertible, then so is its square root $\Gamma _{{\bf X}}^{1/2} $,  as well as the map $f$, whose inverse $f ^{-1} $ is given by $f ^{-1}(\mathbf{x}):=\Gamma _{{\bf X}}^{1/2} \mathbf{x}+ \boldsymbol{\mu}  $. The fact that $f $ is a system isomorphism between $(F, h  )$ and $(\widetilde{F}, \widetilde{h}) $ is a consequence of the equalities \eqref{isomorphic state map}-\eqref{isomorphic readout map}. Parts {\bf (i)} and {\bf (iii)} of Proposition \ref{properties of morphisms solutions} guarantee that if ${\bf X}  $ is the state process associated to $F$ and input ${\bf Z}  $  then so is $\widetilde{{\bf X}} $ defined by $\widetilde{{\bf X}}  _t:=\Gamma_{{\bf X}}^{-1/2} ( \mathbf{X}_t - \boldsymbol{\mu})  $, $t \in \mathbb{Z}_{-} $, with respect to $\widetilde{F }$. The equalities \eqref{mean and covariance new states} immediately follow. \quad $\blacksquare$

\subsection{Proof of Lemma~\ref{rewrite with covs}}
\noindent{
First of all, for any $\tau \in \mathbb{Z}_{-}  $, consider the optimization problems
\begin{equation}
\label{pr1}
(\widehat{{\bf W}}_{{\rm MC}_ \tau}, \hat{a}_{{\rm MC}_ \tau}) =  \stackunder{{\rm arg \ min \ }}{\stackunder{\scriptstyle {\bf W} \in \mathbb{R}^N }{\scriptstyle a \in \mathbb{R} }} {\rm E}\left[\left(\left(T _{-\tau} {\bf Z}\right)_t- {\bf W} ^{\top} U ^F({\bf Z})_t-a\right)^2\right] = \stackunder{{\rm arg \ min \ }}{\stackunder{\scriptstyle {\bf W} \in \mathbb{R}^N }{\scriptstyle a \in \mathbb{R} }} {\rm E}\left[\left(Z_{t+\tau}- {\bf W} ^{\top} {\bf X}_t-a\right)^2\right]
\end{equation}
and 
 \begin{equation}
\label{pr2}
(\widehat{{\bf W}}_{{\rm FC}_ \tau}, \hat{a}_{{\rm FC}_ \tau}) =   \stackunder{{\rm arg \ min \ }}{\stackunder{\scriptstyle {\bf W} \in \mathbb{R}^N }{\scriptstyle a \in \mathbb{R} }} {\rm E}\left[\left( {Z}_t- {\bf W} ^{\top} U ^F(T _{-\tau}({\bf Z}))_t-a\right)^2\right]=\stackunder{{\rm arg \ min \ }}{\stackunder{\scriptstyle {\bf W} \in \mathbb{R}^N }{\scriptstyle a \in \mathbb{R} }} {\rm E}\left[\left( Z_t - {\bf W}^{\top} {\bf X}_{t+\tau} -a\right)^2\right]
\end{equation}
  in  \eqref{definition memory t} and  \eqref{definition forecasting t} of Definition~\ref{memory forecasting capacities}, respectively. It is straightforward to prove by setting equal to zero the derivatives of the objective functions with respect to the optimized parameters (see  Section C in the Technical Supplement of \cite{GHLO2014_capacity} for details) that both \eqref{pr1} and \eqref{pr2}
admit closed-form solutions when $\Gamma_{{\bf X}} $ invertible and they are  given by
\begin{align}
\label{}
\widehat{{\bf W}}_{{\rm MC}_ \tau} = {\rm Cov}( {\bf X}_t,  {\bf X}_t) ^{-1} {\rm Cov}( {\bf X}_t, Z_{t+\tau})= \Gamma_{{\bf X}} ^{-1} {\rm Cov}( {\bf X}_t, Z_{t+\tau}), \quad  \hat{a}_{{\rm MC}_ \tau} = \mu_Z - \widehat{{\bf W}}_{{\rm MC}_ \tau}^\top \boldsymbol{\mu}
\end{align} 
and 
\begin{align}
\label{}
\widehat{{\bf W}}_{{\rm FC}_ \tau} = {\rm Cov}( {\bf X}_{t+\tau},  {\bf X}_{t+\tau}) ^{-1} {\rm Cov}( {\bf X}_{t+\tau}, Z_{t})= \Gamma_{{\bf X}} ^{-1}  {\rm Cov}( {\bf X}_{t+\tau}, Z_{t}), \quad  \hat{a}_{{\rm FC}_ \tau} = \mu_Z - \widehat{{\bf W}}_{{\rm FC}_ \tau}^\top \boldsymbol{\mu}.
\end{align} 
We hence have that 
\begin{align*}
\stackunder{{\rm min}}{\stackunder{\scriptstyle {\bf W} \in \mathbb{R}^N }{\scriptstyle a \in \mathbb{R} }} {\rm E}\left[\left(\left(T _{-\tau} {\bf Z}\right)_t- {\bf W} ^{\top} U ^F({\bf Z})_t-a\right)^2\right] &= {\rm E}\left[\left(Z_{t+\tau} -\widehat{{\bf W}}_{{\rm MC}_ \tau}  ^{\top} {\bf X}_{t}-\hat{a}_{{\rm MC}_ \tau} \right)^2\right] \\
&={\rm Var}(Z_{t+\tau}) -  {\rm Cov}( {\bf X}_{t}, Z_{t+\tau})^\top\Gamma_{{\bf X}} ^{-1} {\rm Cov}( {\bf X}_{t}, Z_{t+\tau}) 
\end{align*}
and
\begin{align*}
 \stackunder{{\rm min \ }}{\stackunder{\scriptstyle {\bf W} \in \mathbb{R}^N }{\scriptstyle a \in \mathbb{R} }} {\rm E}\left[\left( {Z}_t- {\bf W} ^{\top} U ^F(T _{-\tau}({\bf Z}))_t-a\right)^2\right]&= {\rm E}\left[\left(Z_{t} -\widehat{{\bf W}}_{{\rm FC}_ \tau}  ^{\top} {\bf X}_{t+\tau}-\hat{a}_{{\rm FC}_ \tau} \right)^2\right] \\
&={\rm Var}(Z_{t}) -  {\rm Cov}( {\bf X}_{t+\tau}, Z_{t})^\top\Gamma_{{\bf X}} ^{-1} {\rm Cov}( {\bf X}_{t+\tau}, Z_{t}) ,
\end{align*}
which substituted in \eqref{definition memory t} and \eqref{definition forecasting t} and using the variance stationarity of ${\bf Z}$ yield \eqref{capacity formulas when invertible}, as required. \quad $\blacksquare$}

\subsection{Proof of Lemma~\ref{capacity invariant linear morphism}}
\noindent First of all, note that the echo state property hypothesis on  $F _2  $  and part {\bf (ii)} of Proposition \ref{properties of morphisms solutions} imply that  the filter $U^{F _2}$ is well-defined and, moreover, for any ${\bf W} \in \mathbb{R}^{N _2}$ so is $U^{F _1}_{{\bf W} \circ f}$ and 
\begin{equation}
\label{equality after f}
U^{F _1}_{{\bf W} \circ f}=U^{F _2}_{{\bf W}}.
\end{equation}
Let $\langle\cdot , \cdot \rangle_{\mathbb{R}^{N _2}} $ be the Euclidean inner product in $\mathbb{R}^{N _2} $, let $\langle\cdot , \cdot \rangle_{\mathbb{R}^{N _1}}$ be the inner product induced in $\mathbb{R}^{N _1} $ by $\langle\cdot , \cdot \rangle_{\mathbb{R}^{N _2}} $ and the injective map $f$ using \eqref{inner product injective case}, and let $f ^\ast  $ be the corresponding dual map. It is easy to see that the equality \eqref{equality after f} can be rewritten using $f ^\ast  $ as
\begin{equation}
\label{equality after f star}
U^{F _1}_{f ^\ast ({\bf W})}=U^{F _2}_{{\bf W}}, \quad \mbox{for any} \quad {\bf W} \in \mathbb{R}^{N _2}.
\end{equation}
Let now $\tau\in \mathbb{Z}_{-}  $ and let ${\rm MC} _\tau $ be the $\tau$-lag memory capacity of $F_2$ with respect to ${\bf Z} $. By definition and \eqref{equality after f star}
\begin{multline*}
{\rm MC} _\tau = 1- \frac{1}{\text{\rm Var} \left(Z _t\right)}  \stackunder{{\rm min}}{\stackunder{\scriptstyle {\bf W} \in \mathbb{R}^{N_2} }{\scriptstyle a \in \mathbb{R} }} {\rm E}\left[\left(\left(T _{-\tau} {\bf Z}\right)_t-  U ^{F_2}_{{\bf W} }({\bf Z})_t-a\right)^2\right]\\
= 1- \frac{1}{\text{\rm Var} \left(Z _t\right)}  \stackunder{{\rm min}}{\stackunder{\scriptstyle {\bf W} \in \mathbb{R}^{N_2}}{\scriptstyle a \in \mathbb{R} }} {\rm E}\left[\left(\left(T _{-\tau} {\bf Z}\right)_t-  U ^{F_1}_{f ^\ast ({\bf W}) }({\bf Z})_t-a\right)^2\right]\\=
1- \frac{1}{\text{\rm Var} \left(Z _t\right)}  \stackunder{{\rm min}}{\stackunder{\scriptstyle \overline{{\bf W}} \in \mathbb{R}^{N_1}}{\scriptstyle a \in \mathbb{R} }} {\rm E}\left[\left(\left(T _{-\tau} {\bf Z}\right)_t-  U ^{F_1}_{\overline{{\bf W}} }({\bf Z})_t-a\right)^2\right],
\end{multline*}
which coincides with the $\tau$-lag memory capacity of $F_1$. Notice that in the last equality we used the surjectivity of $f ^\ast $ which is a consequence of the injectivity of $f$ (see \eqref{fstarf}). A similar statement can be written for the forecasting capacities. \quad $\blacksquare$

\subsection{Proof of Theorem~\ref{capacity bounds}}
\noindent In the proof we use the following elementary definition and linear algebraic fact: let $\left(V, \langle\cdot , \cdot \rangle_V\right) $ and $\left(W, \langle\cdot , \cdot \rangle_W\right) $ be two inner product spaces and let $f :V \longrightarrow W $ be a linear map between them. The dual map $f ^\ast :W \longrightarrow V $ of $f$  is defined by 
\begin{equation*}
\label{definition dual f}
\langle f ^\ast(\mathbf{w}) , \mathbf{v} \rangle_V=\langle  \mathbf{w} , f(\mathbf{v}) \rangle_W, \quad \mbox{for any} \quad \mathbf{v} \in V, \mathbf{w} \in W.
\end{equation*}
If the map $f$  is injective, then the inner product $\langle\cdot , \cdot \rangle_W $ in $W$ induces an inner product $\langle\cdot , \cdot \rangle_f $ in $V$  via the equality
\begin{equation}
\label{inner product injective case}
\langle\mathbf{v} _1 , \mathbf{v} _2 \rangle_f:=\langle f(\mathbf{v} _1) , f(\mathbf{v} _2) \rangle_W.
\end{equation}
It is easy to see that, in that case, the dual map with respect to the inner product \eqref{inner product injective case} in $V$ and $\langle\cdot , \cdot \rangle_W $ in $W$,  satisfies that 
\begin{equation}
\label{fstarf}
f ^\ast \circ f = \mathbb{I}_V,
\end{equation}
and, in particular, $f ^\ast :W \longrightarrow V $ is surjective.

We now proceed with the proof of the theorem. First of all, since by hypothesis $\Gamma_{{\bf X}}$ is non-singular, Proposition \ref{std rc io} and the second part of Proposition~\ref{properties of morphisms solutions} allow us to replace  the system \eqref{rc state eq}-\eqref{rc readout eq} in  the definitions \eqref{definition memory t} and \eqref{definition forecasting t} by its standardized counterpart whose states $\widetilde{{\bf X}} _t  $ are such that ${\rm E}[\widetilde{\mathbf{X}}_t] = {\bf 0}$ and $ {\rm Cov}( \widetilde{\mathbf{X}}_t, \widetilde{\mathbf{X}}_t)= {\rm E}\left[\widetilde{\mathbf{X}}_t \widetilde{\mathbf{X}}_t ^{\top}\right] = \mathbb{I}_N$. More specifically, if we denote by $\sigma ^2:= \text{\rm Var} \left(Z _t\right)= \gamma(0)$, we can write
\begin{equation*}
{\rm MC} _\tau := 1- \frac{1}{\sigma ^2}  \stackunder{{\rm min}}{\stackunder{\scriptstyle {\bf W} \in \mathbb{R}^N }{\scriptstyle a \in \mathbb{R} }} {\rm E}\left[\left( {Z}_{t+ \tau}- {\bf W} ^{\top} {\bf X}_t-a\right)^2\right]=1- \frac{1}{\sigma ^2}  \stackunder{{\rm min}}{\stackunder{\scriptstyle \widetilde{{\bf W}} \in \mathbb{R}^N }{\scriptstyle \widetilde{a} \in \mathbb{R} }} {\rm E}\left[\left( {Z}_{t+ \tau}- \widetilde{{\bf W}} ^{\top} \widetilde{{\bf X}}_t-\widetilde{a}\right)^2\right],
\end{equation*}
which, using Lemma \ref{rewrite with covs} and the fact that $\Gamma_{\widetilde{{\bf X}}}= \mathbb{I}_N $, can be rewritten as
\begin{equation}
\label{alternative for mc}
{\rm MC} _\tau=\frac{1}{\sigma ^2}{\rm Cov} \left(Z_{t+ \tau}, \widetilde{{\bf X}} _t\right)  {\rm Cov} \left(\widetilde{{\bf X}} _t, Z_{t+ \tau}\right)= \frac{1}{\sigma ^2}\sum_{i=1}^N {\rm E}\left[\widetilde{X} ^i _t (Z_{t+ \tau}- {\rm E}\left[Z_{t+ \tau}\right])\right] ^2.
\end{equation}
Analogously, it is easy to show that 
\begin{equation}
\label{alternative for fc}
{\rm FC} _\tau= \frac{1}{\sigma ^2}\sum_{i=1}^N {\rm E}\left[\widetilde{X} ^i _{t+ \tau} (Z_{t}- {\rm E}\left[Z_{t}\right])\right] ^2.
\end{equation}
If we now define $\widetilde{Z }_t:= (Z _t- {\rm E}\left[Z _t\right])/ \sigma \in L ^2(\Omega, \mathbb{R} )$, for any $t \in \mathbb{Z}_{-} $, it is clear that $\left\|\widetilde{Z }_t\right\|^2_{L ^2}= {\rm E}\left[\widetilde{Z }_t ^2\right]=1$. Moreover, the relation $ {\rm Cov}( \widetilde{\mathbf{X}}_t, \widetilde{\mathbf{X}}_t)= \mathbb{I}_N$ implies that the components $\widetilde{X} _t ^i \in L ^2(\Omega, \mathbb{R} )$, $i \in \left\{1, \ldots, N\right\}$, and that they form an orthonormal set, that is, $\left<\widetilde{X} _t ^i, \widetilde{X} _t ^j\right>_{L ^2}= \delta _{ij} $, where $\delta _{ij} $ stands for Kronecker's delta. { It is actually the properties of the orthogonal projections onto the vector space generated by this orthonormal set that constitute the main technical tool in the proof and that will provide us with the capacity bounds that we are after. }We now separately prove the three parts of the theorem.

\medskip

\noindent\textbf{(i)\ \ } The fact that ${\rm MC} _\tau, {\rm FC} _\tau\geq 0 $ is obvious from \eqref{alternative for mc} and \eqref{alternative for fc}. Let now $\widetilde{S}_t= {\rm span} \left\{\widetilde{X} _t ^1, \ldots, \widetilde{X} _t ^N\right\} \subset L ^2(\Omega, \mathbb{R} )$ and let $\mathbb{P}_{\widetilde{S}_t}:L ^2(\Omega, \mathbb{R} ) \longrightarrow \widetilde{S }_t $ be the corresponding orthogonal projection. Then,
\begin{equation*}
1=\left\|\widetilde{Z }_{t+ \tau}\right\|^2_{L ^2}\geq \left\|\mathbb{P}_{\widetilde{S}_t} \left(\widetilde{Z }_{t+ \tau}\right)\right\|^2_{L ^2}= \left\|\sum_{i=1}^N \widetilde{X} _t ^i\left \langle \widetilde{X} _t ^i,  \widetilde{Z }_{t+ \tau} \right\rangle_{L ^2} \right\|^2_{L ^2}=\sum_{i=1}^N {\rm E}\left[\widetilde{X} ^i _{t} \widetilde{Z}_{t+ \tau}\right] ^2= {\rm MC} _\tau.
\end{equation*}
The inequality ${\rm FC} _\tau\leq 1 $ can be established analogously by considering the projection onto $\widetilde{S}_{t+ \tau  } $ of the vector $\widetilde{Z }_{t} $.

\medskip

\noindent {\bf (ii)} Define first, for any $L \in \mathbb{Z}_{-} $, the vector ${\bf Z} ^L:=({Z} _0- {\rm E}\left[{Z} _0\right], {Z} _{-1}- {\rm E}\left[{Z} _{-1}\right], \ldots, {Z}  _L- {\rm E}\left[{Z} _{L}\right] ) $, and 
\begin{equation}
\label{mc and fc with tildes}
{\rm MC} ^L:=\sum _{\tau=0} ^L {\rm MC} _\tau= \sum _{\tau=0} ^L\sum_{i=1}^N {\rm E}\left[\widetilde{X} ^i _0 \widetilde{Z}_{\tau}\right] ^2, \quad
{\rm FC} ^L:=\sum _{\tau=-1} ^L {\rm FC} _\tau= \sum _{\tau=-1} ^L\sum_{i=1}^N {\rm E}\left[\widetilde{X} ^i _L \widetilde{Z}_{L-\tau}\right] ^2.
\end{equation}
In these equalities we used \eqref{alternative for mc} and \eqref{alternative for fc} as well as the stationarity hypothesis. Now, the properties of the autocovariance function of a second-order stationary process guarantee that the matrix $H ^L  $ is positive semidefinite (see \cite[Theorem 1.5.1]{BrocDavisYellowBook} and \cite{Mukherjee1988} for other properties) and since by hypothesis it is additionally invertible, we can associate to it a square root matrix $\left(H ^L\right)^{1/2} $ that is also invertible. Hence, we define the random vector $\widehat{ {\bf Z}} ^L:=  \left(H ^L\right)^{-1/2} {\bf Z} ^L $, whose components form an orthonormal set in $L ^2(\Omega, \mathbb{R})$. Indeed, for any $i,j \in \left\{1, \ldots, -L+1\right\}$,
\begin{multline}
\label{the Zs are orthogonal}
\langle \widehat{Z} _i^L, \widehat{Z} _j ^L\rangle_{L ^2}= {\rm E}\left[\widehat{Z} _i^L\widehat{Z} _j ^L\right]=
\sum_{k,l=1}^{-L+1} \left(H ^L\right)_{ik}^{-1/2}\left(H ^L\right)_{jl}^{-1/2}{\rm E}\left[({Z} _{-k+1}- {\rm E}\left[{Z} _{-k+1}\right])( {Z} _{-l+1} - {\rm E}\left[{Z} _{-l+1} \right])\right]\\
=
\sum_{k,l=1}^{-L+1} \left(H ^L\right)_{ik}^{-1/2}\left(H ^L\right)_{jl}^{-1/2}\gamma(|k-l|)= \left(\left(H ^L\right)^{-1/2}H ^L \left(H ^L\right)^{-1/2}\right)_{ij}= \delta_{ij}.
\end{multline}
Let  $\widehat{S}_L= {\rm span} \left\{\widehat{Z} _1 ^L, \ldots, \widehat{Z} _{-L+1} ^L\right\} \subset L ^2(\Omega, \mathbb{R} )$ and let $\mathbb{P}_{\widehat{S}_L}:L ^2(\Omega, \mathbb{R} ) \longrightarrow \widehat{S}_L $ be the corresponding orthogonal projection. By \eqref{mc and fc with tildes} and using the definitions introduced earlier we have that 
\begin{align}
{\rm MC} ^L&= \sum _{\tau=0} ^L\sum_{i=1}^N {\rm E}\left[\widetilde{X} ^i _0 \widetilde{Z}_{\tau}\right] ^2={
 \sum _{\tau=0} ^L\sum_{i=1}^N {\rm E}\left[\widetilde{X} ^i _0 \dfrac{Z_{\tau} -{\rm E}\left[ Z_{\tau} \right]}{\sigma}\right] ^2= \frac{1}{\gamma(0)} \sum _{\tau=0} ^L\sum_{i=1}^N {\rm E}\left[\widetilde{X} ^i _0 ({\bf Z}^L)_{-\tau+1}\right] ^2}\nonumber \\
&=\frac{1}{\gamma(0)} \sum _{\tau=0} ^L\sum_{i=1}^N {\rm E}\left[\widetilde{X} ^i _0 \left(\left(H ^L\right)^{1/2} \widehat{{\bf Z}} ^L\right)_{- \tau+1}\right] ^2=
\frac{1}{\gamma(0)} \sum_{i=1}^N \left\|
\left(H ^L\right)^{1/2} {\rm E}\left[\widetilde{X} ^i_0\widehat{{\bf Z}} ^L\right]
\right\|^2.\label{with operator}
\end{align}
Analogously,
\begin{multline}
\label{with operator bis}
{\rm FC} ^L=  \sum _{\tau=-1} ^L\sum_{i=1}^N {\rm E}\left[\widetilde{X} ^i _L \widetilde{Z}_{L- \tau}\right] ^2=
\frac{1}{\gamma(0)} \sum _{\tau=-1} ^L\sum_{i=1}^N {\rm E}\left[\widetilde{X} ^i _L \left(\left(H ^L\right)^{1/2} \widehat{{\bf Z}} ^L\right)_{- L+\tau+1}\right] ^2\\
\leq 
\frac{1}{\gamma(0)} \sum_{i=1}^N \left\|
\left(H ^L\right)^{1/2} {\rm E}\left[\widetilde{X} ^i_L\widehat{{\bf Z}} ^L\right]
\right\|^2.
\end{multline}
Now, by \eqref{the Zs are orthogonal} and \eqref{with operator} we can write that
\begin{multline}
\label{bounding by projection}
{\rm MC} ^L\leq \frac{1}{\gamma(0)} \sum_{i=1}^N 
\vertiii{\left(H ^L\right)^{1/2}}^2 \left\|{\rm E}\left[\widetilde{X} ^i_0\widehat{{\bf Z}} ^L\right]
\right\|^2= 
\frac{1}{\gamma(0)} \sum_{i=1}^N 
\vertiii{\left(H ^L\right)^{1/2}}^2 \left\|\mathbb{P}_{ \widehat{S} _L} \left(\widetilde{X} ^i_0\right)
\right\|_{L ^2}^2\\\leq 
\frac{1}{\gamma(0)} \sum_{i=1}^N 
\vertiii{\left(H ^L\right)^{1/2}} ^2\left\|\widetilde{X} ^i_0
\right\|_{L ^2}^2\leq  \frac{N}{\gamma(0)}|\lambda _{{\rm max}}(H ^L)|=\frac{N}{\gamma(0)} \rho(H ^L).
\end{multline}
An identical inequality can be shown for ${\rm FC} ^L $ using \eqref{with operator bis}, which proves the first inequality in \eqref{bounds memory cap}. The second inequality can be obtained by bounding $\rho(H ^L )$ using Gershgorin's Disks Theorem (see \cite[Theorem 6.1.1 and Corollary 6.1.5]{horn:matrix:analysis}). Indeed, due to this result:
\begin{equation}
\label{first bound of rhoH}
\rho(H ^L )\leq \gamma(0)+\max_{i \in \left\{1, \ldots, -L+1\right\}} 
\left\{
\sum_{\stackon{\scriptstyle j\neq i}{\scriptstyle j\in \left\{1, \ldots, -L+1\right\}}}
\left| H _{ij} ^L\right|
\right\}\leq \gamma (0)+ 2 \sum_{i=1}^{-L}|\gamma(i)|.
\end{equation}

\medskip
\noindent {\bf (iii)} The inequalities in \eqref{inequality summable version} are a consequence of considering $H$ as the infinite symmetric Toeplitz matrix associated to the bi-infinite sequence of autocovariances $\left\{\gamma(j)\right\}_{j \in \mathbb{Z}}$ of ${\bf Z}$. First, when the autocovariance function is absolutely summable then \eqref{spectral density formula} determines the spectral density of ${\bf Z} $ by \cite[Corollary 4.3.2]{BrocDavisYellowBook}. Second, by  \cite[Lemma 6, page 194]{CIT-006}, the spectrum of $H$ is bounded above by the maximum of the function $2 \pi f $, which, using \eqref{bounds memory cap} implies that ${\rm C}\leq \frac{2 \pi N}{\gamma (0)}M _f $. The last inequality is a consequence of 
\begin{equation*}
|f(\lambda) |\leq \frac{1}{2 \pi} \sum _{j=- \infty}^{\infty}|\gamma (j)|= \frac{1}{2 \pi}\left(\gamma(0)+ 2\sum _{j=1}^{\infty}|\gamma (j)|\right)< +\infty, \quad \mbox{for any $\lambda \in [- \pi, \pi]$.\quad $\blacksquare$} 
\end{equation*}

\subsection{Proof of Corollary~\ref{jaeger bound}}
\noindent The first inequality is a straightforward  consequence of \eqref{bounds memory cap} and the fact that for independent inputs $ \gamma(h)=0 $, for all $h \neq 0 $. The second one can be easily obtained from \eqref{alternative for fc}. Indeed, by the causality and the time-invariance \cite[Proposition 2.1]{RC7} of any filter induced by a state-space system of the type \eqref{rc state eq}-\eqref{rc readout eq}, for any $\tau\leq -1 $, the random variables $\widetilde{X} ^i _{t+ \tau} $ and $Z_{t} $ in \eqref{alternative for fc} are independent, and hence
\begin{equation*}
{\rm FC} _\tau= \frac{1}{\sigma ^2}\sum_{i=1}^N {\rm E}\left[\widetilde{X} ^i _{t+ \tau} (Z_{t}- {\rm E}\left[Z_{t}\right])\right] ^2=\frac{1}{\sigma ^2}\sum_{i=1}^N {\rm E}\left[\widetilde{X} ^i _{t+ \tau}\right]^2 {\rm E}\left[ Z_{t}- {\rm E}\left[Z_{t}\right]\right] ^2=0. \quad \blacksquare
\end{equation*}

\subsection{Proof of Proposition~\ref{the linear case with bs}}
\noindent First of all, {recall that the matrix norm $\vertiii{\cdot } $ induced by the Euclidean  norm $\left\|\cdot \right\| $ in $\mathbb{R} ^N$ is defined as $$\vertiii{A} = \stackunder{\rm sup}{{\bf x} \in \mathbb{R} ^N, {\bf x}\neq {\bf 0}} \left\{ \dfrac{\|A{\bf x}\|}{\|{\bf x}\|} \right\} = \sigma_{{\rm max}}(A).$$ 
It follows from this definition   that the condition $\sigma_{{\rm max}}(A)<1 $ implies that the state map $F(\mathbf{x}, z)$  in \eqref{definition linear system} is a contraction on the first entry. Indeed, for any $\mathbf{x} _1, \mathbf{x} _2 \in \mathbb{R}^N $, $z \in \mathbb{R} $, we have
\begin{equation*}
\left\|F(\mathbf{x}_1, z)-F(\mathbf{x}_2, z)\right\|= \left\|A(\mathbf{x}_1-\mathbf{x}_2)\right\|\leq \vertiii{A} \left\|\mathbf{x}_1-\mathbf{x}_2\right\|<\left\|\mathbf{x}_1-\mathbf{x}_2\right\|.
\end{equation*}
} Second, as the input process takes values on a compact set, there exists a compact subset $D_N\subset \mathbb{R}^N  $ (see \cite[Remark 2]{RC10}) such that the restriction $F: D_N \times D \longrightarrow D _N $ satisfies the hypotheses of Proposition \ref{ESP for reservoir maps with compact target} and of Corollary \ref{stationary implies stationary}. This implies that the system associated to $F$ has the echo state property, as well as the stationarity of the filter output ${\bf X} =U^{A, \mathbf{C}}({\bf Z})$ and of the joint processes in the statement. We recall that, in this case,
\begin{equation*}
{\bf X}_t =U^{A, \mathbf{C}}({\bf Z})_t=\sum _{j=0} ^{\infty}A ^j\mathbf{C} {Z}_{t-j}, \quad \mbox{for any} \quad t \in \mathbb{Z}_{-}.
\end{equation*}
We now show that the output process is also square-integrable and hence covariance stationary. Indeed, let ${\bf X} _t ^n:= \sum _{j=0} ^{n}A ^j\mathbf{C} {Z}_{t-j} $, $n \in \mathbb{N} $. Given that by hypothesis $D  $ is compact, there exists $M> 0 $ such that $D \subset [-M,M] $ and hence
\begin{equation*}
\left\|{\bf X} _t^n\right\|= \left\| \sum _{j=0} ^{n}A ^j\mathbf{C} {Z}_{t-j}\right\|\leq
M\left\|\mathbf{C}\right\| \sum _{j=0} ^{n}\vertiii{A} ^j\leq \frac{M\left\|\mathbf{C}\right\|}{1- \sigma_{{\rm max}}(A)}.
\end{equation*}
Now the Bounded Convergence Theorem guarantees that
\begin{equation*}
\left\|{\bf X} _t\right\|_{L ^2}={\rm E}\left[\left\|{\bf X} _t\right\|^2\right]^{\frac{1}{2}}= \lim_{n \rightarrow \infty} {\rm E}\left[\left\|{\bf X} _t^n\right\|^2\right]^{\frac{1}{2}}\leq \frac{M\left\|\mathbf{C}\right\|}{1- \sigma_{{\rm max}}(A)}< \infty.
\end{equation*}
The Cauchy-Schwarz inequality implies that the components of $\Gamma_{{\bf X}}= {\rm E}\left[{\bf X}_t {\bf X}_t^{\top} \right]$ are also finite and hence, using the notation introduced in Proposition \ref{std rc io} and the invertibility hypothesis on $\Gamma_{{\bf X}}$, the standarized states $\widetilde{{\bf X}} _t $ are given by
\begin{equation}
\label{standarized state 1}
\widetilde{{\bf X}} _t=\Gamma_{{\bf X}}^{-1/2}{\bf X}_t=\Gamma_{{\bf X}}^{-1/2}\sum _{j=0} ^{\infty}A ^j\mathbf{C} {Z}_{t-j}.
\end{equation}
Additionally, when the autocovariance function of the input is absolutely summable, then the spectral density $f$ of ${\bf Z} $ defined in \eqref{spectral density formula} belongs to the so-called Wiener class and, moreover, if the hypothesis on it in the statement is satisfied, then  the two matrices $H$ (semi-infinite) and $\overline{H} $ (doubly infinite) are invertible (see \cite[Theorem 11]{CIT-006}).

\medskip

\noindent {\bf (i)}  Let $\overline{{\bf Z}}:= \left({ Z} _0, { Z} _1, \ldots\right) $  and let $\widehat{{\bf Z}}=H^{-1/2}\overline{{\bf Z}}$. An argument similar to \eqref{the Zs are orthogonal} shows that
$\langle \widehat{Z} _i, \widehat{Z} _j \rangle_{L ^2}= \delta_{ij}$. Moreover, \eqref{standarized state 1} implies that
\begin{equation*}
\widetilde{{\bf X}} _0=\Gamma_{{\bf X}}^{-1/2}\sum _{k=0} ^{\infty}A ^k \mathbf{C} {Z}_{-k}=\Gamma_{{\bf X}}^{-1/2}\sum _{k=0} ^{\infty}A ^k \mathbf{C} \left( H^{1/2}\widehat{{\bf Z}}\right)_{k+1}=\Gamma_{{\bf X}}^{-1/2}\sum _{k=0} ^{\infty}\sum_{j=1} ^{\infty}A ^k \mathbf{C}  H^{1/2}_{k+1,j}\widehat{{Z}}_j,
\end{equation*}
which, using the definition in \eqref{definition b1} can be rewritten componentwise as
\begin{equation}
\label{XtildeinZhatbasis}
\widetilde{{X}}_0^i=\sum_{j=1}^{\infty} {B}_i^j \widehat{{Z}}_j, \qquad i \in \left\{1, \ldots, N\right\}.
\end{equation}
The scalars ${B}_i^j $ are defined in \eqref{definition b1} and by \eqref{XtildeinZhatbasis} coincide with the (unique) coefficients that determine the expansion of $\widetilde{{X}}_0^i $ on the orthonormal basis $\left\{\widehat{Z} _t\right\}_{t \in \mathbb{Z}_{-}}$. This implies, in particular, that ${B}_i^j = \langle \widetilde{{X}}_0^i, \widehat{{Z}}_j\rangle_{L ^2}$ and, by the Cauchy-Schwarz inequality, that all these coefficients are finite.

We now show that the $L ^2 $-orthonormality of the components $\widetilde{{X}}_0^i $ of $\widetilde{{\bf X}} _0$ implies the $ \ell ^2 $-orthonormality of the vectors $\left\{\mathbf{B} _1, \ldots, \mathbf{B} _N \right\}  \in \ell^2_+(\mathbb{R}) $ whose components we just showed are finite. Indeed, for any $i,j \in \left\{1, \ldots, N\right\} $, the equality \eqref{XtildeinZhatbasis}  and the Parseval identity imply that
\begin{equation}
\label{ortonormality in l2}
\delta_{ij}=\langle  \widetilde{{X}}_0^i  ,  \widetilde{{X}}_0^j  \rangle _{L^2}
=\sum_{k=1}^{\infty} {B}_i^k   {B}_j^k= \langle \mathbf{B}_i, \mathbf{B}_j\rangle _{\ell^2_+}.
\end{equation}
We now prove \eqref{expression memory linear case}. First, taking the limit $L \rightarrow \infty $ in \eqref{with operator} we write,
\begin{multline*}
{\rm MC} = 
\frac{1}{\gamma(0)} \sum_{i=1}^N \left\|
H ^{1/2} {\rm E}\left[\widetilde{X} ^i_0\widehat{{\bf Z}} \right]
\right\|^2=
\frac{1}{\gamma(0)} \sum_{i=1}^N \left\|
H ^{1/2} \sum_{j=1}^{\infty}B _i^j{\rm E}\left[\widehat{{Z}}_j \widehat{{\bf Z}}\right]\right\|^2\\
=
\frac{1}{\gamma(0)} \sum_{i=1}^N \left\|
H ^{1/2} \mathbf{B} _i\right\|^2_{\ell^2}= \frac{1}{\gamma(0)}\sum_{i=1}^N \langle \mathbf{B} _i, H \mathbf{B}_i \rangle _{\ell^2}.
\end{multline*}

\medskip

\noindent {\bf (ii)}  Let now $\overline{{\bf Z}}:= \left(\ldots, Z _{-1}, { Z} _0, { Z} _1, \ldots\right) $  and let $\widehat{{\bf Z}}=\overline{H}^{-1/2}\overline{{\bf Z}}$. An argument similar to \eqref{the Zs are orthogonal} shows that
$\langle \widehat{Z} _i, \widehat{Z} _j \rangle_{L ^2}= \delta_{ij}$ for any $i,j \in \Bbb Z$. Also, in this case, \eqref{standarized state 1} implies that
\begin{equation*}
\widetilde{{\bf X}} _0=\Gamma_{{\bf X}}^{-1/2}\sum _{k=0} ^{\infty}\sum_{j=- \infty} ^{\infty}A ^k \mathbf{C}  \overline{H}^{1/2}_{-k,j}\widehat{{Z}}_j,
\end{equation*}
which, using the definition in \eqref{definition b2} can be rewritten componentwise as
\begin{equation*}
\widetilde{{X}}_0^i=\sum_{j=- \infty}^{\infty} {B}_i^j \widehat{{Z}}_j, \qquad i \in \left\{1, \ldots, N\right\}.
\end{equation*}
As in \eqref{ortonormality in l2}, 
the $L ^2 $-orthonormality of the components $\widetilde{{X}}_0^i $ of $\widetilde{{\bf X}} _0$ implies the $ \ell ^2 $-orthonormality of the vectors $\left\{\mathbf{B} _1, \ldots, \mathbf{B} _N \right\}  \in \ell^2(\mathbb{R}) $. Finally, in order to establish \eqref{expression forecasting linear case}, we use the stationarity hypothesis to rewrite the expression of the forecasting capacity in \eqref{alternative for fc} as 
\begin{multline*}
{\rm FC} = 
\frac{1}{\gamma(0)} \sum_{\tau=1}^{\infty}\sum_{i=1}^N 
 {\rm E}\left[\widetilde{X} ^i_0{Z} _\tau\right]
^2=\frac{1}{\gamma(0)} \sum_{\tau=1}^{\infty}\sum_{i=1}^N 
 {\rm E}\left[\sum_{j=- \infty}^{\infty} {B}_i^j \widehat{{Z}}_j\left(\overline{H}^{1/2}\widehat{{\bf Z}}\right)_\tau\right]
^2
\\
=\frac{1}{\gamma(0)} \sum_{\tau=1}^{\infty}\sum_{i=1}^N 
\left(
\sum_{j,l=- \infty}^{\infty}B_i ^j \overline{H}^{1/2}_{\tau,l} {\rm E}\left[\widehat{{Z}}_j \widehat{{Z}}_l\right]
\right)^2=
\frac{1}{\gamma(0)} \sum_{\tau=1}^{\infty}\sum_{i=1}^N 
\left(
\sum_{j=- \infty}^{\infty}B_i ^j \overline{H}^{1/2}_{\tau,j} 
\right)^2\\
=
\frac{1}{\gamma(0)} \sum_{\tau=1}^{\infty}\sum_{i=1}^N 
\left(
 \overline{H}^{1/2} \mathbf{B} _i 
\right)_\tau^2= \frac{1}{\gamma(0)}\sum_{i=1}^N 
\left\|\mathbb{P}_{\mathbb{Z}^+}\left(\overline{H}^{1/2}\mathbf{B} _i\right)\right\|^2_{\ell^2}. \quad \blacksquare
\end{multline*}

\subsection{Proof of Corollary~\ref{jaeger bound for linear}}
\noindent The first equality in \eqref{traditional capacities linear} is a straightforward consequence of \eqref{expression memory linear case} and of the fact that for white noise inputs $H= \gamma(0) \mathbb{I}_{\ell^2_+(\mathbb{R})}$. The equality ${\rm FC}=0 $ follows from the fact that $\overline{H}= \gamma(0) \mathbb{I}_{\ell^2(\mathbb{R})}$ and that  $B ^j_i=0 $, for any $i\in \left\{1, \ldots, N\right\} $ and any $j \in \Bbb Z^+ $, by \eqref{definition b2}. Then, by \eqref{expression forecasting linear case},
\begin{equation*}
{\rm FC}= \frac{1}{\gamma(0)}\sum_{i=1}^N 
\left\|\mathbb{P}_{\mathbb{Z}^+}\left(\overline{H}^{1/2}\mathbf{B} _i\right)\right\|^2_{\ell^2}=\sum_{i=1}^N 
\left\|\mathbb{P}_{\mathbb{Z}^+}\left(\mathbf{B} _i\right)\right\|^2_{\ell^2}=0. \qquad \blacksquare
\end{equation*}

\subsection{Proof of Proposition~\ref{capacities and kalman}}
\noindent \textbf{(i)} Using the notation introduced in the statement notice that:
\begin{equation}
\label{kalman vectors in eigenvalue basis}
A \mathbf{C}=\sum_{i=1}^N c _i\lambda _i\mathbf{v} _i, \, 
A^2 \mathbf{C}=\sum_{i=1}^N c _i\lambda _i^2\mathbf{v} _i, \, \ldots ,\,
A^N \mathbf{C}=\sum_{i=1}^N c _i\lambda _i^N\mathbf{v} _i.
\end{equation}
Since by hypothesis $\vertiii{A}_2=\sigma_{{\rm max}}(A)<1 $, the spectral radius $\rho(A) $ of $A$ satisfies that $\rho(A)\leq \sigma_{{\rm max}}(A)<1 $, and hence:
\begin{equation}
\label{BmatrixDef} 
B:=\sum_{j=0}^{\infty}A ^j \mathbf{C}\mathbf{C} ^{\top} \left(A ^j\right)^{\top}=\sum_{k=0}^{\infty}\sum_{i,j=1}^N \lambda _i^k \lambda_j^k c _i c  _j \mathbf{v} _i \mathbf{v} _j^{\top}=\sum_{i,j=1}^N \frac{c _ic _j}{1- \lambda _i\lambda _j}\mathbf{v} _i \mathbf{v} _j^{\top},
\end{equation}
which shows that in the matrix basis $\left\{\mathbf{v} _i \mathbf{v} _j^{\top}\right\}_{i,j \in \left\{1, \ldots, N\right\}} $, the matrix $B$ has components $\overline{B} _{ij}:= \frac{c _ic _j}{1- \lambda _i\lambda _j} $ or, equivalently,
\begin{equation}
\label{definition b bar}
\overline{B}:= \overline{\mathbf{C}}\overline{\mathbf{C}} ^{\top}\odot D, \quad \mbox{with} \quad \overline{\mathbf{C}}= (c _1, \ldots, c _N)^{\top}, \quad \mbox{and $D$ defined by} \quad D_{ij}:=\frac{1}{1- \lambda _i\lambda _j},
\end{equation}
for any $ i,j \in \left\{1, \ldots,N\right\}$; the symbol $\odot $ stands for componentwise matrix multiplication (Hadamard product). Let $P\in \mathbb{M} _N$ be the invertible change-of-basis matrix between $\{ \mathbf{v}_1, \ldots, \mathbf{v}_N\}  $ and the canonical basis $\left\{\mathbf{e}_1, \ldots, \mathbf{e}_N\right\}$, that is, for any $i \in \left\{1, \ldots,N\right\} $, we have
that $\mathbf{v} _i=\sum_{k=1}^NP _{ik}\mathbf{e}_k $. It is easy to see that $B=P ^{\top} \overline{B} P$ and hence the invertibility of $B$  (and hence of $\Gamma_{{\bf X}} $) is equivalent to the invertibility of $\overline{B} $ in \eqref{definition b bar}, which we now characterize. 

In order to provide an alternative expression for $\overline{B} $, recall that for any two vectors $\mathbf{v}, \mathbf{w} \in \mathbb{R}^N $, we can write $\mathbf{v}\odot \mathbf{w}= {\rm diag}(\mathbf{v}) \mathbf{w} $, where ${\rm diag}(\mathbf{v}) \in \mathbb{M}_N $  is the diagonal matrix that has the entries of the vector $\mathbf{v}  $ in the diagonal. Using this fact and the Hadamard product trace property (see \cite[Lemma 5.1.4, page 305]{Horn:Johnson}) we have that
\begin{multline*}
\left\langle \mathbf{v}, \overline{B} \mathbf{w} \right\rangle=\left\langle \mathbf{v}, (\overline{\mathbf{C}}\overline{\mathbf{C}} ^{\top}\odot D) \mathbf{w}\right\rangle={\rm trace} \left(\mathbf{v}^{\top} (\overline{\mathbf{C}}\overline{\mathbf{C}} ^{\top}\odot D) \mathbf{w}\right)=
{\rm trace} \left( (\overline{\mathbf{C}}\overline{\mathbf{C}} ^{\top}\odot D) \mathbf{w}\mathbf{v}^{\top}\right)\\
={\rm trace} \left( (\overline{\mathbf{C}}\overline{\mathbf{C}} ^{\top}\odot \mathbf{v} \mathbf{w}^{\top})D ^{\top} \right)=
{\rm trace} \left( (\overline{\mathbf{C}}\odot \mathbf{v})( \overline{\mathbf{C}} \odot \mathbf{w})^{\top}D ^{\top} \right)\\
={\rm trace} \left({\rm diag}(\overline{\mathbf{C}})\mathbf{v} \mathbf{w}^{\top}{\rm diag}(\overline{\mathbf{C}})D ^{\top}
\right)=
\left\langle \mathbf{v},  {\rm diag}(\overline{\mathbf{C}}) D{\rm diag}(\overline{\mathbf{C}})\mathbf{w} \right\rangle.
\end{multline*}
Since $\mathbf{v}, \mathbf{w} \in \mathbb{R}^N $ are arbitrary, this equality allows us to conclude that $\overline{B}={\rm diag}(\overline{\mathbf{C}}) D{\rm diag}(\overline{\mathbf{C}}) $ and hence $\overline{B} $ is invertible if and only if both ${\rm diag}(\overline{\mathbf{C}}) $  and  $D $ are. The regularity of ${\rm diag}(\overline{\mathbf{C}}) $ is equivalent to requiring that all the entries of the vector $\overline{\mathbf{C}} $ are non-zero. Regarding $D $, it can be shown by induction on the matrix dimension $N $, that
\begin{equation*}
{\rm det}(D)= \left(-1\right)^N \frac{\prod_{i<j=1}^N \left(\lambda_i- \lambda _j\right)^2}{\prod_{i<j=1}^N \left(\lambda _i\lambda _j-1\right)^2\prod_{i=1}^N(\lambda_i^2-1)}.
\end{equation*}
Consequently, $D$ is invertible if and only if ${\rm det}(D)\neq 0 $, which is equivalent to all the elements in the spectrum $\sigma(A)$ being distinct.

\medskip

\noindent\textbf{(ii)} The condition on the vectors $\left\{A \mathbf{C}, A ^2\mathbf{C}, \ldots, A ^N \mathbf{C}\right\}$ forming a basis of $ \mathbb{R}^N$ is equivalent to the invertibility of the  matrix $ \widehat{R(A, \mathbf{C})}:= \left(A \mathbf{C}| A ^2\mathbf{C}| \cdots| A ^N \mathbf{C}\right)$. It is easy to see using \eqref{kalman vectors in eigenvalue basis} that 
\begin{equation}
\label{equivalence with hat}
\widehat{R(A, \mathbf{C})}= P ^{\top}\overline{R(A, \mathbf{C})},
\end{equation}
where $P$ is the invertible change-of-basis matrix in the previous point and $\overline{R(A, \mathbf{C})} $ is given by
\begin{equation}
\label{matrix m bar}
\overline{R(A, \mathbf{C})}:=
\left(
\begin{array}{cccc}
c _1\lambda_1&c _1\lambda_1 ^2 &\cdots &c _1\lambda_1 ^N\\
c _2\lambda_2&c _2\lambda_2^2 &\cdots &c _2\lambda_2^N\\
\vdots & \vdots &\ddots &\vdots\\
c _N\lambda_N&c _N\lambda_N^2 &\cdots &c _N\lambda_N^N
\end{array}
\right).
\end{equation}
Indeed, for any $i, j \in \left\{1, \ldots, N\right\} $,
\begin{multline*}
 \widehat{R(A, \mathbf{C})} _{ij}= \left(A ^j \mathbf{C}\right)_i= \left(\sum_{k=1}^N c _k A ^j \mathbf{v}_k\right)_i=\left(\sum_{k=1}^N c _k \lambda_k ^j \mathbf{v}_k\right)_i=
\left(\sum_{k,l=1}^N c _k \lambda_k ^jP_{kl} \mathbf{e}_l\right)_i\\=
\left(\sum_{k,l=1}^N \overline{R(A, \mathbf{C})}_{kj} P_{kl} \mathbf{e}_l\right)_i=
\left(\sum_{l=1}^N \left(P ^{\top}\overline{R(A, \mathbf{C})}\right)_{lj} \mathbf{e}_l\right)_i=
\sum_{l=1}^N \left(P ^{\top}\overline{R(A, \mathbf{C})}\right)_{lj} \mathbf{e}_i ^{\top}\mathbf{e}_l=
\left(P ^{\top}\overline{R(A, \mathbf{C})}\right)_{ij}, 
\end{multline*}
which proves \eqref{equivalence with hat}.
Now, using induction on the matrix dimension $N$, it can be shown that 
\begin{equation*}
{\rm det}(\overline{R(A, \mathbf{C})})=\prod_{i=1}^N c _i \lambda_i\prod_{i<j=1}^N(\lambda_i- \lambda _j).
\end{equation*}
The invertibility of $\overline{R(A, \mathbf{C})} $ (or, equivalently, the invertibility of $\widehat{R(A, \mathbf{C})}$) is equivalent to all the coefficients $c _i $ and all the eigenvalues $\lambda _i $ being non-zero (so that $\prod_{i=1}^N c _i \lambda_i $ is non-zero) and all the elements in $\sigma(A)$ being distinct (so that $\prod_{i<j=1}^N(\lambda_i- \lambda _j) $ is non-zero).

\medskip

\noindent\textbf{(iii)} The Kalman controllability condition on the vectors
$\left\{\mathbf{C}, A \mathbf{C}, \ldots, A ^{N-1} \mathbf{C}\right\}$ forming a basis of $ \mathbb{R}^N$ is equivalent to the invertibility of the {\bfi  controllability} or {\bfi  reachability} matrix $R(A, \mathbf{C}):= \left(\mathbf{C}| A \mathbf{C}| \cdots| A ^{N-1} \mathbf{C}\right)$ (see \cite{sontag:book} for this terminology). Following the same strategy that we used to prove \eqref{equivalence with hat}, it is easy to see that $R(A, \mathbf{C})= P ^{\top}\widetilde{R(A, \mathbf{C})} $, where
\begin{equation*}
\widetilde{R(A, \mathbf{C})}:=
\left(
\begin{array}{cccc}
c _1&c _1 \lambda _1 &\cdots &c _1 \lambda _1^{N-1}\\
c _2&c _2 \lambda _2 &\cdots &c _2 \lambda _2^{N-1}\\
\vdots & \vdots &\ddots &\vdots\\
c _N&c _N \lambda _N &\cdots &c _N \lambda _N^{N-1}
\end{array}
\right).
\end{equation*}
The matrix $\widetilde{R(A, \mathbf{C})} $ has the same rank as $\overline{R(A, \mathbf{C})} $ since it can be obtained from $\widetilde{R(A, \mathbf{C})} $ via elementary matrix operations, namely, by dividing each row $i$ of $\overline{R(A, \mathbf{C}) } $ by the corresponding eigenvalue $\lambda _i  $ which is by hypothesis non-zero. \quad $\blacksquare$
{
\begin{remark} 
\normalfont
The matrix $B$ in \eqref{BmatrixDef}, which is instrumental for the proof of part {\bf (i)} of the proposition, is related to the objects introduced in \cite{Tino2019}, in particular to the positive semi-definite symmetric matrix  $Q$  corresponding  to the temporal kernel associated to the linear dynamical system and defined as $Q_{i,j} = \mathbf{C}^\top (A^{j-1})^\top A^{i-1} \mathbf{C}$, $i,j\in\{1,\ldots, \tau\}$, $\tau \in \mathbb{N} $. More specifically,  if one defines $B_{\tau}$ as the covariance matrix of the states process associated to the truncated solution of the linear state system \eqref{definition linear system} with  white noise as inputs, that is, $B_{\tau}:=\sum_{j=0}^{\tau} \mathbf{C} ^{\top}  \left(A ^j\right)^{\top}A ^j \mathbf{C}$, $\tau \in \mathbb{N} $, then it is easy to see that ${\rm trace}(B_{\tau}) = \sum^{\tau}_{j = 1} \mathbf{C}^\top (A^j)^\top A^j \mathbf{C}= \sum^{\tau}_{j = 1}Q_{j+1,j+1} = {\rm trace}(Q^d)$ with $Q^d \in \mathbb{S}_{\tau} $ the diagonal matrix with the same elements on the main diagonal as $Q$. The results in \cite{Tino2019} provide bounds for the elements $Q_{i,j}$,  $i,j\in\{1,\ldots, \tau\}$, whenever the state map is constructed with a randomly generated connectivity matrix $A$ and input matrix ${\mathbf C}$, for different choices of architectures and distributions. Obviously, the analysis of the diagonal entries of $Q$ for those situations is valid for the diagonal elements of $B_{\tau}$ and hence, since $\Gamma_{\mathbf X} = \gamma(0) B$, for large $\tau$ also illustrates the behavior of the variances of the states process of  linear state systems.
\end{remark}}

\subsection{Proof of Theorem~\ref{expression capacity linear system}}
\noindent\textbf{(i)} We first show that $\ker R(A, \mathbf{C}) ^{\top} \subset \ker \Gamma_{\mathbf{X}} $. Let $\mathbf{v} \in \ker R(A, \mathbf{C}) ^{\top} $. This implies that 
\begin{equation}
\label{powers zero 1}
\mathbf{C}^{\top}(A ^j)^{\top} \mathbf{v}= {\bf 0},  \  \mbox{for all $j \in \left\{0, \ldots, N-1\right\} $.}
\end{equation}
Now, by the Hamilton-Cayley Theorem \cite[Theorem 2.4.3.2]{horn:matrix:analysis}, for any $j\geq N $, there exist constants $\left\{\beta_0^j, \ldots, \beta_{N-1}^j\right\}$ such that $A ^j=\sum _{i=0}^{N-1} \beta _i^j A^i $ which, together with \eqref{powers zero 1}, implies that equality holds for all $j \in \mathbb{N}  $. Recall now that by Proposition \ref{capacities and kalman} part {\bf (i)}, $\Gamma_{{\bf X}}= \gamma(0)\sum_{j=0}^{\infty}A ^j \mathbf{C}\mathbf{C} ^{\top} \left(A ^j\right)^{\top}$, and hence we can conclude that $\Gamma_{{\bf X}}(\mathbf{v})=\gamma(0)\sum_{j=0}^{\infty}A ^j \mathbf{C}\mathbf{C} ^{\top} \left(A ^j\right)^{\top} \mathbf{v}= {\bf 0} $, that is, $\mathbf{v} \in \ker \Gamma_{\mathbf{X}} $. Conversely, if $\mathbf{v} \in \ker \Gamma_{\mathbf{X}} $, we have that $0=\langle\mathbf{v}, \Gamma_{\mathbf{X}}(\mathbf{v})\rangle= \gamma (0)\sum_{j=0}^{\infty}  \|\mathbf{C} ^{\top} \left(A ^j\right)^{\top} \mathbf{v} \|^2$, which implies that $\mathbf{C} ^{\top} \left(A ^j\right)^{\top} \mathbf{v}= {\bf 0} $, necessarily, for any $j \in \mathbb{N}  $ and hence $\mathbf{v} \in \ker R(A, \mathbf{C}) ^{\top} $.

\medskip

\noindent {\bf (ii)} The system associated to $\overline{F} $ has the echo state property because for any $\overline{\mathbf{x}} \in X $,
\begin{equation*}
\left\|\overline{A} \overline{\mathbf{x}}\right\|^2= \left\|\pi _XA i_X(\overline{\mathbf{x}})\right\|^2\leq \left\|A i_X(\overline{\mathbf{x}})\right\|^2,
\end{equation*}
which implies that $\vertiii{\overline{A}}_2\leq \vertiii{ {A}}_2=\sigma_{{\rm max}}(A)<1$.

\medskip

\noindent {\bf (iii)} We first show that for any $\overline{\mathbf{x}} \in X $ and $z \in D  $ we have that $A i _X \overline{\mathbf{x}}+ \mathbf{C} z  \in X $. Indeed, as $\overline{ \mathbf{x}} \in X $, there exist constants $\left\{\alpha_0, \ldots, \alpha_{N-1}\right\}$ such that $\overline{\mathbf{x}} =\sum_ {i=0}^{N-1} \alpha_i A^i \mathbf{C} $ and hence 
\begin{multline}
\label{x is stable}
A i _X \overline{\mathbf{x}}+ \mathbf{C} z=\sum_ {i=0}^{N-1} \alpha_i A^{i+1} \mathbf{C} + \mathbf{C}z=\mathbf{C}z+\sum_ {i=1}^{N-1} \alpha_{i-1} A^{i} \mathbf{C}+ \alpha_{N-1}A ^N\mathbf{C}\\=
\mathbf{C}z+\sum_ {i=1}^{N-1} \alpha_{i-1} A^{i} \mathbf{C}+ \alpha_{N-1}\sum_{j=0}^{N-1}\beta_j^NA ^j \mathbf{C} \in X,
\end{multline}
where the constants $\left\{\beta ^N _0, \ldots , \beta ^N _{N-1}\right\}$  satisfy that $A ^N=\sum _{i=0}^{N-1} \beta _i^N A^i $ and, as above, are a byproduct of the Hamilton-Cayley Theorem \cite[Theorem 2.4.3.2]{horn:matrix:analysis}. We now show that $i _X $  is a system equivariant map between $\overline{F}  $  and $F $.  For any $\overline{\mathbf{x}} \in X $ and $z \in D  $,
\begin{equation*}
i _X \left(\overline{F}(\overline{\mathbf{x}}, z)\right)= i _X \pi_X \left(A i _X(\overline{\mathbf{x}})+ \mathbf{C} z \right)=A i _X(\overline{\mathbf{x}})+ \mathbf{C} z=F(i _X(\overline{\mathbf{x}}),z),
\end{equation*}
where the second equality holds because $i _X \circ \pi _X|_{X}= \mathbb{I}_X  $  and, by \eqref{x is stable}, $A i _X \overline{\mathbf{x}}+ \mathbf{C} z \in X $. 

\medskip

\noindent {\bf (iv)} We  prove the identity ${\rm rank}\, R(\overline{A}, \overline{\mathbf{C}})=r $ by showing that 
\begin{equation*}
i _X \left( {\rm span} \left\{\overline{\mathbf{C}}, \overline{A}\overline{\mathbf{C}}, \ldots, \overline{A}^{r-1}\overline{\mathbf{C}}\right\}\right)=X.
\end{equation*}
First, using that $i _X \circ \pi _X|_{X}= \mathbb{I}_X  $ and that by the Hamilton-Cayley Theorem $A ^j\mathbf{C} \in X $, for all $j \in \mathbb{N} $, it is easy to conclude that for any $\mathbf{v}=\sum_{i=1}^r \alpha _i \overline{A}^{i-1}\overline{\mathbf{C}} $  we have that
\begin{equation*}
i _X(\mathbf{v})= i _X \left(\sum_{i=1}^r \alpha _i \overline{A}^{i-1}\overline{\mathbf{C}}\right)
=\sum_{i=1}^r \alpha _i  {A}^{i-1} {\mathbf{C}},
\end{equation*}
which shows that $i _X \left( {\rm span} \left\{\overline{\mathbf{C}}, \overline{A}\overline{\mathbf{C}}, \ldots, \overline{A}^{r-1}\overline{\mathbf{C}}\right\}\right) \subset X$. Conversely, let $\mathbf{v}=\sum_{i=1}^{N}\alpha _iA^{i-1}\mathbf{C} \in X$. Using again that $i _X \circ \pi _X|_{X}= \mathbb{I}_X  $, we can write that,
\begin{multline*}
\mathbf{v}=\sum_{i=1}^{N}\alpha _iA^{i-1}\mathbf{C}=\sum_{i=1}^{N}\alpha_i\underbrace{\left(i _X \pi _XA\right)\cdots \left(i _X \pi _XA\right)}_\text{$i-1 $-times} \left(i _X \pi _X\mathbf{C}\right)=
i _X \left(\sum_{i=1}^N \alpha _i \overline{A}^{i-1}\overline{\mathbf{C}}\right)\\
=
i _X \left(\sum_{i=1}^{r} \alpha _i \overline{A}^{i-1}\overline{\mathbf{C}}+
\sum_{j=r+1}^{N}\sum_{k _j=1}^r \alpha _j \beta ^j_{l} \overline{A}^{l-1}\overline{\mathbf{C}}
\right),
\end{multline*}
which clearly belongs to $i _X \left( {\rm span} \left\{\overline{\mathbf{C}}, \overline{A}\overline{\mathbf{C}}, \ldots, \overline{A}^{r-1}\overline{\mathbf{C}}\right\}\right)$. The constants $\beta ^j_{l} $ are obtained, again, by using the Hamilton-Cayley Theorem. Finally, the invertibility of $ \Gamma_{\overline{{\bf X}}} $ is a consequence of Proposition \ref{capacities and kalman} and the hypothesis that the elements of the spectrum $\sigma(\overline{A}) $ are non-zero.

\medskip

\noindent {\bf (v)} The statement in part {\bf (iii)} that we just proved and Lemma \ref{capacity invariant linear morphism} imply that the memory ${\rm MC} $ and forecasting ${\rm FC}$ capacities of $F$  with respect to ${\bf Z} $ coincide with those of $\overline{F} $ with respect to ${\bf Z} $. Now, the statement in part {\bf (iv)}  and Corollary \ref{jaeger bound for linear} imply that those capacities coincide with $r$ and $0$, respectively. Finally, as $r = {\rm rank}\,R(A, \mathbf{C})$, the claim follows. \quad $\blacksquare$

\medskip

\noindent {\bf Acknowledgments:} We thank the editor and two anonymous referees for their careful reading and interesting comments that have improved the paper. L Gonon and JPO acknowledge partial financial support  coming from the Research Commission of the Universit\"at Sankt Gallen and the Swiss National Science Foundation (grant number 200021\_175801/1).  JPO acknowledges partial financial support of the French ANR ``BIPHOPROC" project (ANR-14-OHRI-0002-02). The authors thank the hospitality and the generosity of the FIM at ETH Zurich and the Division of Mathematical Sciences of the Nanyang Technological University, Singapore, where a significant portion of the results in this paper were obtained. 

\noindent
\addcontentsline{toc}{section}{Bibliography}
\bibliographystyle{wmaainf}

\end{document}